\documentclass[11pt]{article}

\usepackage[in]{fullpage}
\usepackage{times}
\usepackage{helvet}
\usepackage{courier}
\usepackage{graphicx}
\usepackage[cmex10]{amsmath}
\usepackage[mathscr]{euscript}
\usepackage{bm}
\usepackage{amsfonts}
\usepackage{tikz}
\usepackage[mathscr]{euscript}

\usepackage{hyperref}

\newtheorem{theorem}{Theorem}[section]
\newtheorem{lemma}[theorem]{Lemma}

\newtheorem{definition}[theorem]{Definition}

\begin{document}

\title{Estimating Structured Vector Autoregressive Models}

\author{Igor Melnyk \\
Dept of Computer Science \& Engineering\\
University of Minnesota, Twin Cities\\
melnyk@cs.umn.edu
\and
Arindam Banerjee\\
Dept of Computer Science \& Engineering\\
University of Minnesota, Twin Cities\\
banerjee@cs.umn.edu}



\maketitle

\begin{abstract}
While considerable advances have been made in estimating high-dimensional structured models from independent data using Lasso-type models, limited progress has been made for settings when the samples are dependent. We consider estimating structured VAR (vector auto-regressive models), where the structure can be captured by any suitable norm, e.g., Lasso, group Lasso, order weighted Lasso, sparse group Lasso, etc. In VAR setting with correlated noise, although there is strong dependence over time and covariates, we establish bounds on the non-asymptotic estimation error of structured VAR parameters. Surprisingly, the estimation error is of the same order as that of the corresponding Lasso-type estimator with independent samples, and the analysis holds for any norm. Our analysis relies on results in generic chaining, sub-exponential martingales, and spectral representation of VAR models. Experimental results on synthetic data with a variety of structures as well as real aviation data are presented, validating theoretical results.
\end{abstract}

\section{Introduction}
\label{sec:intro}
The past decade has seen considerable progress on approaches to structured parameter estimation, especially in the linear regression setting, where one considers regularized estimation problems of the form:
\begin{align}
\label{eq:basic}
\hat{\bm{\beta}} = \underset{\bm{\beta} \in \mathbb{R}^q}{\text{argmin}}~\frac{1}{M} \| \mathbf{y} - Z \bm{\beta} \|^2_2 + \lambda_M R(\bm{\beta})~,
\end{align}
where $\{(y_i,z_i), i=1,\ldots,M\}$, $y_i \in \mathbb{R}, z_i \in \mathbb{R}^q$, such that $\mathbf{y} = [y_1^T,\ldots,y_M^T]^T$ and $Z=[z_1^T,\ldots, z_M^T]^T$, is the training set of $M$ independently and identically distributed (i.i.d.) samples, $\lambda_M > 0$ is a regularization parameter, and $R(\cdot)$ denotes a suitable norm \cite{tibshirani96, zou05, yuan06}. Specific choices of $R(\cdot)$ lead to certain types of structured parameters to be estimated. For example, the decomposable norm $R(\bm{\beta}) = \| \bm{\beta} \|_1$ yields Lasso, estimating sparse parameters, $R(\bm{\beta}) = \|\bm{\beta}\|_G$ gives Group Lasso, estimating group sparse parameters, and $R(\bm{\beta}) = \|\bm{\beta}\|_{owl}$, the ordered weighted $L_1$ norm (OWL) \cite{bbsc13}, gives sorted $L_1$-penalized estimator, clustering correlated regression parameters \cite{fino12}. Non-decomposable norms, such as $K$-support norm \cite{arfs12} or overlapping group sparsity norm \cite{jacob09} can be used to uncover more complicated model structures. Theoretical analysis of such models, including sample complexity and non-asymptotic bounds on the estimation error rely on the design matrix $Z$, usually assumed (sub)-Gaussian with independent rows, and the specific norm $R(\cdot)$ under consideration \cite{raskutti10,rudelson13}.  Recent work has generalized such estimators to work with any norm \cite{negahbann12,banerjee14} with i.i.d.~rows in $Z$.

The focus of the current paper is on structured estimation in vector auto-regressive (VAR) models \cite{lutkepohl07}, arguably the most widely used family of multivariate time series models. VAR models have been applied widely, ranging from describing the behavior of economic and financial time series \cite{tsay05} to modeling the dynamical systems  \cite{ljung98} and estimating brain function connectivity \cite{valdes05}, among others.  A VAR model of order $d$ is defined as
\begin{align}
\label{eq:vard}
x_t = A_1 x_{t-1} + A_2 x_{t-2} + \cdots + A_d x_{t-d} + \epsilon_t~,
\end{align}
where $x_t \in \mathbb{R}^p$ denotes a multivariate time series, $A_k \in \mathbb{R}^{p \times p}, k=1,\ldots,d$ are the parameters of the model, and $d \geq 1$ is the order of the model. In this work, we assume that the noise $\epsilon_t \in \mathbb{R}^p$ follows a Gaussian distribution, $\epsilon_t \sim \mathcal{N}(0,\Sigma)$, with $\mathbb{E}(\epsilon_t \epsilon_t^T ) = \Sigma$ and $\mathbb{E}(\epsilon_t \epsilon_{t+\tau}^T) = 0$, for $\tau \neq 0$. The VAR process is assumed to be stable and stationary \cite{lutkepohl07}, while the noise covariance matrix $\Sigma$ is assumed to be positive definite with bounded largest eigenvalue, i.e., $\Lambda_{\text{min}}(\Sigma) > 0$ and $\Lambda_{\text{max}}(\Sigma) < \infty$.

In the current context, the parameters $\{A_k\}$ are assumed to be structured, in the sense of having low values according to a suitable norm $R(\cdot)$. We consider a general setting where {\em any} norm can be applied to the rows $A_k(i,:) \in \mathbb{R}^p$ of $A_k$, allowing the possibility of different norms being applies to different rows of $A_k$, and different norms for different parameter matrices $A_{k}, k=1,\ldots,d$. Choosing $L_1$-norm $\| A_k(i,:) \|_1$ for all rows and all parameter matrices is a simple special case of our setting. We discuss certain other choices in Section \ref{sec:regEstimator}, and discuss related results in Section~\ref{sec:results}. In order to estimate the parameters, one can consider regularized estimators of the form \eqref{eq:basic}, where $y_i$ and $z_i$ correspond to $x_t$ in the VAR setting. Unfortunately, unlike $(y_i, z_i)$ in \eqref{eq:basic}, the $x_t$ are {\em far from independent}, having strong dependence across time and correlated across dimensions. As a result, existing results
from the rich literature on regularized estimators for structured problems \cite{zhyu06, wain09, meyu09} cannot be directly applied to get sample complexities and estimation error bounds in VAR models.

The rest of the paper is organized as follows. In Section \ref{sec:related} we review the related work on structured VAR estimation. In Section \ref{sec:model} we present the estimation problem for the VAR model and in Section \ref{sec:theory} we establish the main results of our analysis on the VAR estimation guarantees. We present experimental results in Section \ref{sec:results} and conclude in Section \ref{sec:conc}. The proofs and other detailed discussions can be found in Appendices A, B, C and D.

\section{Related Work}
\label{sec:related}
In recent literature, the problem of estimating structured VAR models has been considered for the special case of $L_1$ norm. \cite{han13} analyzed a constrained estimator based on the Dantzig selector \cite{candes07}, and established the recovery results for the special case of $L_1$ norm. \cite{song11} considered a regularized VAR estimation problem under Lasso and Group Lasso penalties and derived oracle inequalities for the prediction error and estimation accuracy. However, their analysis is for the case when the dimensionality of the problem is fixed with respect to the sample size. Moreover, they employed an assumption on the dependency structure in the VAR, thus limiting the sample correlation issues mentioned earlier.

The work of \cite{koca15}  studied regularized Lasso-based estimator while allowing for problem dimensionality to grow with sample size, utilizing suitable martingale concentration inequalities to analyze dependency structure. \cite{loh11} considered $L_1$ VAR estimation for first order models ($d=1$) assuming $\|A_1\|_2 < 1$, and the analysis was not extended to the general case of $d>1$. In recent work, \cite{bami15} considered a VAR Lasso estimator and established the sample complexity and error bounds by building on the prior work of \cite{loh11}. Their approach exploits the spectral properties of a general VAR model of order $d$, providing insights on the dependency structure of the VAR process. However, in line with the existing literature, the analysis was tailored to the special case of $L_1$ norm, thus limiting its generality.

Compared to the existing literature, our results are substantially more general since the results and analysis apply to {\em any} norm $R(\cdot)$. One may wonder---given the popularity of $L_1$ norm, why worry about other norms? Over the past decade, considerable effort has been devoted to generalize $L_1$ norm based results to other norms \cite{negahbann12, csbc12, banerjee14, fino12}. Our work obviates the need for a similar exercise for VAR models. Further, some of these norms have found key niche in specific application areas e.g., \cite{zlny12, ywsh15}. From a technical perspective, one may also wonder---once we have the result for $L_1$ norm, why should not the extension to other norms be straightforward? A key technical aspect of the estimation error analysis boils down to getting sharp concentration bounds for $R^*(Z^T \epsilon)$, where $R^*(\cdot)$ is the dual norm of $R(\cdot)$, $Z$ is the design matrix, and $\epsilon$ is the noise \cite{banerjee14}. For the special case of $L_1$, the dual norm is $L_{\infty}$, and one can use {\em union bound} to get the required concentration. In fact, this is exactly how the analysis in \cite{bami15} was done. For general norms, the union bound is inapplicable. Our analysis is based on a considerably more powerful tool, {\em generic chaining} \cite{talagrand06}, yielding an analysis applicable to any norm, and producing results in terms of geometric properties, such as Gaussian widths \cite{leta13}, of sets related to the norm. Results for specific norms can then be obtained by plugging in suitable bounds on the Gaussian widths \cite{crpw12, chba15}. We illustrate the idea by recovering known bounds for Lasso and Group Lasso, and obtaining new results for Spare Group Lasso and OWL norms. Finally, in terms of the core technical analysis, the application of generic chaining to the VAR estimation setting is not straightforward. In the VAR setting, generic chaining has to consider a stochastic process derived from sub-exponential martingale difference sequence (MDS). We first generalize the classical Azuma-Hoeffding inequality applicable to sub-Gaussian MDSs to get an Azuma-Bernstein inequality for sub-exponential MDSs. Further, we use suitable representations of Talagrand's $\gamma$-functions \cite{talagrand06} in the context of generic chaining to obtain bounds on $R^*(Z^T \epsilon)$ in terms of the Gaussian width $w(\Omega_R)$ of the unit norm ball $\Omega_R =\{ u \in \mathbb{R}^{dp} | R(u) \leq 1 \}$. Our estimation error bounds in the VAR setting are {\em exactly of the same order} as Lasso-type models in the i.i.d.~setting implying, surprisingly, that the strong temporal dependency in the VAR setting has no adverse effect on the estimation.

\section{Structured VAR Model}
\label{sec:model}
In this section we formulate structured VAR estimation problem and discuss its properties, which are essential in characterizing sample complexity and error bounds.


\subsection{Regularized Estimator}
\label{sec:regEstimator}
To estimate the parameters of the VAR model, we transform the model in \eqref{eq:vard} into the form suitable for regularized estimator \eqref{eq:basic}. Let $(x_0, x_1, \ldots, x_T)$ denote the $T+1$ samples generated by the stable VAR model in \eqref{eq:vard}, then stacking them together we obtain
\begin{align*}
\begin{bmatrix}
x_d^T\\
x_{d+1}^T\\
\vdots\\
x_T^T\\
\end{bmatrix}
=
\begin{bmatrix}
x_{d-1}^T & \hspace{-5pt}x_{d-2}^T & \hdots & \hspace{-5pt} x_{0}^T\\
x_{d}^T & \hspace{-5pt}x_{d-1}^T & \hdots & \hspace{-5pt}x_{1}^T\\
\vdots  &\hspace{-5pt} \vdots & \ddots &\hspace{-5pt} \vdots \\
x_{T-1}^T &\hspace{-5pt} x_{T-2}^T &\hdots & \hspace{-5pt}x_{T-d}^T
\end{bmatrix}
\begin{bmatrix}
A_1^T\\
A_2^T\\
\vdots\\
A_d^T
\end{bmatrix}
+
\begin{bmatrix}
\epsilon_d^T\\
\epsilon_{d+1}^T\\
\vdots\\
\epsilon_T^T
\end{bmatrix}
\end{align*}
which can also be compactly written as
\begin{align}
\label{eq:VARmatrix}
Y = XB + E,
\end{align}
where $Y \in \mathbb{R}^{N \times p}$, $X \in \mathbb{R}^{N \times dp}$, $B \in \mathbb{R}^{dp\times p}$, and $E \in \mathbb{R}^{N \times p}$ for $N = T-d+1$. Vectorizing (column-wise) each matrix in \eqref{eq:VARmatrix}, we get
\begin{align*}
\text{vec}(Y) &= (I_{p\times p}\otimes X)\text{vec}(B) + \text{vec}(E)\\
\mathbf{y} &= Z\bm{\beta} + \bm{\epsilon},
\end{align*}
\noindent where $\mathbf{y} \in \mathbb{R}^{Np}$, $Z= (I_{p\times p}\otimes X) \in \mathbb{R}^{Np\times dp^2}$, $\bm{\beta} \in \mathbb{R}^{dp^2}$, $\bm{\epsilon} \in \mathbb{R}^{Np}$, and $\otimes$ is the Kronecker product. The covariance matrix of the noise $\bm{\epsilon}$ is now $\mathbb{E}[\bm{\epsilon}\bm{\epsilon}^T] = \Sigma\otimes I_{N\times N}$. Consequently, the regularized estimator takes the form
\begin{align}
\label{eq:VARestimator}
\hat{\bm{\beta}} = \underset{\bm{\beta} \in \mathbb{R}^{dp^2}}{\text{argmin}} \frac{1}{N}||\mathbf{y} - Z\mathbf{\bm{\beta}}||_2^2 + \lambda_NR(\bm{\beta}),
\end{align}
\noindent where $R(\bm{\beta})$ can be any vector norm, separable along the rows of matrices $A_k$. Specifically, if we denote $\bm{\beta} = [\beta_1^T \ldots \beta_p^T]^T$ and $A_k(i,:)$ as the row of matrix $A_k$ for $k=1,\ldots, d$, then our assumption is equivalent to
\begin{align}
\label{eq:decompNorm}
R(\bm{\beta})\hspace{-2pt} =\hspace{-3pt} \sum_{i=1}^p\hspace{-3pt}R\big(\beta_i\big) \hspace{-3pt}=\hspace{-3pt} \sum_{i=1}^p\hspace{-3pt}R\bigg(\hspace{-2pt}\Big[A_1(i,:)^T \hspace{-3pt}\ldots\hspace{-1pt} A_d(i,:)^T\hspace{-0pt}\Big]^T\hspace{-2pt}\bigg).
\end{align}
To reduce clutter and without loss of generality, we assume the norm $R(\cdot)$ to be the same for each row $i$. 
Since the analysis decouples across rows, it is straightforward to extend our analysis to the case when a different norm is used for each row of $A_k$, e.g., $L_1$ for row one, $L_2$ for row two, $K$-support norm \cite{arfs12} for row three, etc. Observe that within a row, the norm need not be decomposable across columns.

The main difference between the estimation problem in \eqref{eq:basic} and the formulation in \eqref{eq:VARestimator} is the strong dependence between the samples $(x_0, x_1, \ldots, x_T)$, violating the i.i.d.~assumption on the data $\{(y_i,z_i), i=1,\ldots,Np\}$. In particular, this leads to the correlations between the rows and columns of matrix $X$ (and consequently of $Z$). To deal with such dependencies, following \cite{bami15}, we utilize the spectral representation of the autocovariance of VAR models to control the dependencies in matrix $X$. 

\subsection{Stability of VAR Model}
\label{sec:stability}
Since VAR models are (linear) dynamical systems, for the analysis we need to establish conditions under which the VAR model \eqref{eq:vard} is stable, i.e., the time-series process does not diverge over time.  
For understanding stability, it is convenient to rewrite VAR model of order $d$ in \eqref{eq:vard} as an equivalent VAR model of order $1$
\begin{align}
\label{eq:var1}
\begin{bmatrix}
x_t\\
x_{t-1}\\
\vdots\\
x_{\hspace{-1pt}t\hspace{-1pt}-\hspace{-1pt}(\hspace{-1pt}d\hspace{-1pt}-\hspace{-1pt}1\hspace{-1pt})}
\end{bmatrix}
\hspace{-4pt}
 =
\hspace{-4pt}
\underbrace{
\begin{bmatrix}
A_1 & \hspace{-5pt}A_2 & \hspace{-5pt}\ldots & \hspace{-3pt}A_{d-1} &\hspace{-3pt} A_d\\
I & \hspace{-5pt}0 & \hspace{-5pt}\ldots & \hspace{-3pt}0 & \hspace{-3pt}0\\
0 & \hspace{-5pt}I & \hspace{-5pt}\ldots & \hspace{-3pt}0 & \hspace{-3pt}0\\
\vdots& \hspace{-5pt}\vdots & \hspace{-5pt}\ddots & \hspace{-3pt}\vdots & \hspace{-3pt}\vdots\\
0 & \hspace{-5pt}0 & \hspace{-5pt}\ldots & \hspace{-3pt}I & \hspace{-3pt}0
\end{bmatrix}}_{\mathbf{A}}
\hspace{-5pt}
\begin{bmatrix}
x_{t\hspace{-1pt}-\hspace{-1pt}1}\\
x_{t\hspace{-1pt}-\hspace{-1pt}2}\\
\vdots\\
x_{t\hspace{-1pt}-\hspace{-1pt}d}
\end{bmatrix}
\hspace{-4pt}
+
\hspace{-4pt}
\begin{bmatrix}
\epsilon_t\\
0\\
\vdots\\
0
\end{bmatrix}
\end{align}
where $\mathbf{A}\in \mathbb{R}^{dp\times dp}$. 
Therefore, VAR process is stable if all the eigenvalues of $\mathbf{A}$ satisfy $\text{det}(\lambda I_{dp\times dp}- \mathbf{A}) = 0$ for $\lambda \in \mathbb{C}$, $|\lambda| < 1$. Equivalently, if expressed in terms of original parameters $A_k$, stability is satisfied if $\text{det}(I-\sum_{k=1}^dA_k\frac{1}{\lambda^{k}}) = 0$ (see Appendix \ref{sec:stability_append} for more details).

\subsection{Properties of Data Matrix $X$}
\label{sec:Xproperties}

In what follows, we analyze the covariance structure of matrix $X$ in \eqref{eq:VARmatrix} using spectral properties of VAR model (see Appendix \ref{sec:covariance_append} for additional details). The results will then be used in establishing the high probability bounds for the estimation guarantees in problem \eqref{eq:VARestimator}.

Define any row of $X$ as $X_{i,:} \in \mathbb{R}^{dp}$, $1\leq i \leq N$. Since we assumed that $\epsilon_t \sim \mathcal{N}(0,\Sigma)$, it follows that each row is distributed as $X_{i,:} \sim \mathcal{N}(0, C_{\mathsf{X}})$, where the covariance matrix $C_{\mathsf{X}} \in \mathbb{R}^{dp\times dp}$ is the same for all $i$
\begin{align}
\label{eq:covX}
C_{\mathsf{X}} = \begin{bmatrix}
\Gamma(0) & \Gamma(1) & \ldots & \Gamma(d-1)\\
\vspace{-4pt}
\Gamma(1)^T & \Gamma(0) & \ldots & \Gamma(d-2)\\
\vspace{-2pt}
\vdots & \vdots & \ddots & \vdots\\
\Gamma(d-1)^T & \Gamma(d-2)^T & \ldots & \Gamma(0)\\
\end{bmatrix},
\end{align}
where $\Gamma(h) = \mathbb{E}(x_tx_{t+h}^T) \in \mathbb{R}^{p\times p}$. It turns out that since $C_{\mathsf{X}}$ is a block-Toeplitz matrix, its eigenvalues can be bounded as (see \cite{blockToeplitz11})
\begin{align}
\label{eq:eigbounds}
\underset{
\begin{subarray}{c}
\vspace{.2em}
  1 \leq j \leq p\\
  \omega \in [0, 2\pi]
\end{subarray}}{\inf} \Lambda_j[\rho(\omega)] \leq \underset{1 \leq k \leq dp}{\Lambda_k[C_\mathsf{X}]} \leq
\underset{
\begin{subarray}{c}
\vspace{.2em}
  1 \leq j \leq p\\
  \omega \in [0, 2\pi]
\end{subarray}}{\sup} \Lambda_j[\rho(\omega)],
\end{align}
where  $\Lambda_k[\cdot]$ denotes the $k$-th eigenvalue of a matrix and for $i = \sqrt{-1}$, $\rho(\omega) = \sum_{h=-\infty}^{\infty} \Gamma(h) e^{-hi\omega}, \quad \omega \in [0, 2\pi]$,
is the spectral density, i.e., a Fourier transform of the autocovariance matrix $\Gamma(h)$. The advantage of utilizing spectral density is that it has a closed form expression (see Section 9.4 of \cite{priestley81})
\begin{align*}
\rho(\omega) \hspace{-3pt}=\hspace{-3pt} \left(I\hspace{-3pt}-\hspace{-3pt}\sum_{k=1}^dA_ke^{-ki\omega}\hspace{-3pt}\right)^{-1}\hspace{-10pt}\Sigma\left[\left(I\hspace{-3pt}-\hspace{-3pt}\sum_{k=1}^dA_ke^{-ki\omega}\right)^{-1}\hspace{-0pt}\right]^{*},
\end{align*}
where $*$ denotes a Hermitian of a matrix. Therefore, from \eqref{eq:eigbounds} we can establish the following lower bound
\begin{align}
\label{eq:autocovBoundss}
\Lambda_{\min}[C_\mathsf{X}] \geq \Lambda_{\min}(\Sigma)/\Lambda_{\max}(\mathscr{A}) = \mathscr{L} ,
\end{align}
where we defined $\Lambda_{\max}(\mathscr{A}) = \underset{\omega \in [0,2\pi]}{\max}\Lambda_{\max}(\mathscr{A}(\omega))$ for
\begin{align}
\label{eq:caligA}
\mathscr{A}(\omega) \hspace{-2pt}=\hspace{-2pt} \left(I\hspace{-0pt}-\hspace{-0pt}\sum_{k=1}^dA_k^Te^{ki\omega}\right)\left(I\hspace{-0pt}-\hspace{-0pt}\sum_{k=1}^dA_ke^{-ki\omega}\right),
\end{align}
see Appendix \ref{sec:VARcovx} for additional details.

In establishing high probability bounds we will also need information about a vector $q = Xa \in \mathbb{R}^{N}$ for any $a \in \mathbb{R}^{dp}$, $\|a\|_2=1$.  Since each element $X_{i,:}^Ta \sim \mathcal{N}(0, a^TC_{\mathsf{X}}a)$, it follows that $q \sim \mathcal{N}(0, Q_a)$ with a covariance matrix $Q_a \in \mathbb{R}^{N\times N}$. It can be shown (see Appendix \ref{sec:rowDistrr} for more details) that $Q_a$ can be written as
\begin{align}
\label{eq:Qa}
Q_a = (I\otimes a^T)C_\mathcal{U}(I\otimes a),
\end{align}
where $C_\mathcal{U} = \mathbb{E}(\mathcal{U}\mathcal{U}^T)$ for $\mathcal{U} =
\begin{bmatrix}
X_{1,:}^T & \ldots & X_{N,:}^T
\end{bmatrix} ^T\in \mathbb{R}^{Ndp}$
%
%
which is obtained from matrix $X$ by stacking all the rows in a single vector, i.e, $\mathcal{U} = \text{vec}(X^T)$. In order to bound eigenvalues of $C_\mathcal{U}$ (and consequently of $Q_a$), observe that $\mathcal{U}$ can be viewed as a vector obtained  by stacking $N$ outputs from VAR model in \eqref{eq:var1}. Similarly as in \eqref{eq:eigbounds}, if we denote the spectral density of the VAR process in \eqref{eq:var1} as $\rho_{\mathsf{X}}(\omega) = \sum_{h=-\infty}^{\infty} \Gamma_{\mathsf{X}}(h) e^{-hi\omega}$, $\omega \in [0, 2\pi]$, where  $\Gamma_{\mathsf{X}}(h) = \mathbb{E}[X_{j,:}X_{j+h,:}^T] \in \mathbb{R}^{dp\times dp}$, then we can write
\begin{align*}
\underset{
\begin{subarray}{c}
\vspace{.2em}
  1 \leq l \leq dp\\
  \omega \in [0, 2\pi]
\end{subarray}}{\inf} \Lambda_l[\rho_{\mathsf{X}}(\omega)] \leq
\underset{1 \leq k \leq Ndp}{
 \Lambda_k[C_\mathcal{U}]} \leq
\underset{
\begin{subarray}{c}
\vspace{.2em}
  1 \leq l \leq dp\\
  \omega \in [0, 2\pi]
\end{subarray}}{\sup} \Lambda_l[\rho_{\mathsf{X}}(\omega)].
\end{align*}
The closed form expression of spectral density is
\begin{align*}
\rho_{\mathsf{X}}(\omega) = \left(I - \mathbf{A}e^{-i\omega}\right)^{-1}\Sigma_\mathcal{E} \left[\left(I - \mathbf{A}e^{-i\omega}\right)^{-1}\right]^{*},
\end{align*}
where $\Sigma_\mathcal{E}$ is the covariance matrix of a noise vector and $\mathbf{A}$ are as defined in expression \eqref{eq:var1}. Thus, an upper bound on $C_\mathcal{U}$ can be obtained as $\Lambda_{\max}[C_\mathcal{U}] \leq
\frac{\Lambda_{\max}(\Sigma)}{\Lambda_{\min}(\bm{\mathscr{A}})}$, where we defined $\Lambda_{\min}(\bm{\mathscr{A}}) = \underset{\omega \in [0,2\pi]}{\min}\Lambda_{\min}(\bm{\mathscr{A}}(\omega))$ for
\begin{align}
\label{eq:caligAbold}
\bm{\mathscr{A}}(\omega) = \left(I - \mathbf{A}^Te^{i\omega}\right)\left(I - \mathbf{A}e^{-i\omega}\right).
\end{align}
Referring back to covariance matrix $Q_a$ in \eqref{eq:Qa}, we get
\begin{align}
\label{eq:Q}
\Lambda_{\max}[Q_a] \leq
\Lambda_{\max}(\Sigma)/\Lambda_{\min}(\bm{\mathscr{A}}) = \mathscr{M}.
\end{align}

We note that for a general VAR model, there might not exist closed-form expressions for $\Lambda_{\max}(\mathscr{A})$ and $\Lambda_{\min}(\bm{\mathscr{A}})$. However, for some special cases there are results establishing the bounds on these quantities (e.g., see Proposition~2.2 in \cite{bami15}).

\section{Regularized Estimation Guarantees}
\label{sec:theory}
Denote by $\Delta = \hat{\bm{\beta}} - \bm{\beta}^*$ the error between the solution of optimization problem \eqref{eq:VARestimator} and $\bm{\beta}^*$, the true value of the parameter. The focus of our work is to determine conditions under which the optimization problem in \eqref{eq:VARestimator} has guarantees on the accuracy of the obtained solution, i.e., the error term is bounded: $||\Delta||_2 \leq \delta$ for some known $\delta$. To establish such conditions, we utilize the framework of \cite{banerjee14}. Specifically, estimation error analysis is based on the following known results adapted to our settings. The first one characterizes the restricted error set $\Omega_{E}$, where the error $\Delta$ belongs.   
\begin{lemma}
\label{lemmaReg}
Assume that 
\begin{align}
\label{eq:lambdaBound}
\lambda_N \geq rR^*\left[\frac{1}{N}Z^T\bm{\epsilon}\right], 
\end{align}
\noindent for some constant $r > 1$, where $R^*\left[\frac{1}{N}Z^T\bm{\epsilon}\right]$ is a dual form of the vector norm $R(\cdot)$, which is defined as $R^*[\frac{1}{N}Z^T\bm{\epsilon}] = \underset{R(U) \leq 1}{\sup} \left<\frac{1}{N}Z^T\bm{\epsilon}, U\right>$, for $U \in \mathbb{R}^{dp^2}$, where $U = [u_1^T, u_2^T, \ldots, u_p^T]^T$ and $u_i \in \mathbb{R}^{dp}$. Then the error vector $\|\Delta\|_2$ belongs to the set
\begin{align}
\label{eq:errorset}
\Omega_E \hspace{-3pt}=\hspace{-3pt} \left\{ \Delta \in \mathbb{R}^{dp^2} \Big| R(\bm{\beta}^* \hspace{-3pt}+\hspace{-3pt} \Delta) \leq R(\bm{\beta}^*) \hspace{-3pt}+\hspace{-3pt} \frac{1}{r}R(\Delta)\right\}.
\end{align}
\end{lemma}
The second condition in \cite{banerjee14} establishes the upper bound on the estimation error.
\begin{lemma}
\label{lemmaRE}
Assume that the restricted eigenvalue (RE) condition holds 
\begin{align}
\label{eq:reBound}
\frac{||Z\Delta||_2}{||\Delta||_2} \geq \sqrt{\kappa N}, 
\end{align}
for $\Delta \in \textnormal{cone}(\Omega_E)$ and some constant $\kappa > 0$, where $\textnormal{cone}(\Omega_E)$ is a cone of the error set, then
\begin{align}
\label{eq:errorBound}
||\Delta||_2 \leq \frac{1+r}{r}\frac{\lambda_N}{\kappa}\Psi(\textnormal{cone}(\Omega_E)),
\end{align}
\noindent where $\Psi(\textnormal{cone}(\Omega_E))$ is a norm compatibility constant, defined as $\Psi(\textnormal{cone}(\Omega_E)) = \underset{U \in \text{\textnormal{cone}}(\Omega_E)}{\sup}\frac{R(U)}{||U||_2}$.
\end{lemma}
Note that the above error bound is deterministic, i.e., if \eqref{eq:lambdaBound} and \eqref{eq:reBound} hold, then the error satisfies the upper bound in \eqref{eq:errorBound}. However, the results are defined in terms of the quantities, involving $Z$ and $\bm{\epsilon}$, which are random. Therefore, in the following we establish high probability bounds on the regularization parameter in \eqref{eq:lambdaBound} and RE condition in \eqref{eq:reBound}.

\subsection{High Probability Bounds}

In this Section we present the main results of our work, followed by the discussion on their properties and illustrating some special cases based on popular Lasso and Group Lasso regularization norms. In Section \ref{sec:proofScetch} we will present the main ideas of our proof technique, with all the details delegated to the Appendices \ref{sec:lambdaBound_append} and \ref{sec:re_append}. 

To establish lower bound on the regularization parameter $\lambda_N$, we derive an upper bound on $R^*[\frac{1}{N}Z^T\bm{\epsilon}] \leq \alpha$, for some $\alpha>0$, which will establish the required relationship $\lambda_N \geq \alpha \geq R^*[\frac{1}{N}Z^T\bm{\epsilon}]$. 

\begin{theorem}
\label{theoremReg}
Let $\Omega_R = \{u \in \mathbb{R}^{dp}| R(u) \leq 1\}$, and define $w(\Omega_R) = \mathbb{E}[\underset{u\in \Omega_R}{\sup}\left< g, u\right>]$ to be a Gaussian width of set $\Omega_R$ for $g\sim \mathcal{N}(0,I)$. 
 For any $\epsilon_1 > 0$ and $\epsilon_2 >0$ with probability at least $1- c \exp(-\min(\epsilon_2^2, \epsilon_1) + \log(p))$ we can establish that 
\begin{align*}
R^*\left[\frac{1}{N}Z^T\bm{\epsilon}\right]\leq \left(c_2(1\hspace{-3pt}+\hspace{-3pt}\epsilon_2)\frac{w(\Omega_R)}{\sqrt{N}} + c_1(1\hspace{-3pt}+\hspace{-3pt}\epsilon_1)\frac{w^2(\Omega_R)}{N^2}\right)
\end{align*}
where $c$, $c_1$ and $c_2$ are positive constants.
\end{theorem}
To establish restricted eigenvalue condition, we will show that $\underset{\Delta \in \textnormal{cone}(\Omega_E)}{\inf}\frac{||(I_{p\times p}\otimes X)\Delta||_2}{||\Delta||_2} \geq \nu$, for some $\nu > 0$ and then set $\sqrt{\kappa N} = \nu$.
\begin{theorem}
\label{theoremRE}
Let $\Theta = \textnormal{cone}(\Omega_{E_j}) \cap S^{dp-1}$, where $S^{dp-1}$ is a unit sphere. The error set $\Omega_{E_j}$ is defined as $\Omega_{E_j} = \left\{ \Delta_j \in \mathbb{R}^{dp} \Big| R(\beta_j^* + \Delta_j) \leq R(\beta_j^*) + \frac{1}{r}R(\Delta_j)\right\}$, for $r > 1$, $j=1,\ldots, p$, and $\Delta = [\Delta_1^T, \ldots, \Delta_p^T]^T$, for $\Delta_j$ is of size $dp\times 1$, and $\bm{\beta}^* = [\beta_1^{*T} \ldots \beta_p^{*T}]^T$, for $\beta_j^* \in \mathbb{R}^{dp}$. The set $\Omega_{E_j}$ is a part of the decomposition in $\Omega_E = \Omega_{E_1} \times \cdots \times \Omega_{E_p}$ due to the assumption on the row-wise separability of norm $R(\cdot)$ in \eqref{eq:decompNorm}. Also define $w(\Theta) = \mathbb{E}[\underset{u\in \Theta}{\sup}\left< g, u\right>]$ to be a Gaussian width of set $\Theta$ for $g\sim \mathcal{N}(0,I)$ and $u \in \mathbb{R}^{dp}$. Then with probability at least $1 - c_1\exp(-c_2\eta^2 + \log(p))$, for any $\eta > 0$
\begin{align*}
\underset{\Delta \in \textnormal{cone}(\Omega_E)}{\inf}\hspace{-0pt}\frac{||(I_{p\times p}\hspace{0pt}\otimes\hspace{-0pt} X)\Delta||_2}{||\Delta||_2} \geq \nu,
\end{align*}
where $\nu = 
\sqrt{N\mathscr{L}} - 2\sqrt{\mathscr{M}} - c w(\Theta) - \eta$ and $c$, $c_1$, $c_2$ are positive constants, and $\mathscr{L}$ and $\mathscr{M}$ are defined in \eqref{eq:autocovBoundss} and \eqref{eq:Q}.
\end{theorem}

\subsection{Discussion}
\label{sec:discussion}
From Theorem \ref{theoremRE}, we can choose $\eta = \frac{1}{2}\sqrt{N\mathscr{L}}$ and set $\sqrt{\kappa N} = \sqrt{N\mathscr{L}} - 2\sqrt{\mathscr{M}} - c w(\Theta) - \eta$ and since $\sqrt{\kappa N} > 0$ must be satisfied, we can establish a lower bound on the number of samples $N$
\begin{align}
\label{eq:Nbound}
\sqrt{N} > \frac{2\sqrt{\mathscr{M}} + c w(\Theta)}{\sqrt{\mathscr{L}}/2} = \mathcal{O}(w(\Theta)).
\end{align}
Examining this bound and using \eqref{eq:autocovBoundss} and \eqref{eq:Q}, we can conclude that the number of samples needed to satisfy the restricted eigenvalue condition is smaller if $\Lambda_{\min}(\bm{\mathscr A})$ and $\Lambda_{\min}(\Sigma)$ are larger and $\Lambda_{\max}({\mathscr A})$ and $\Lambda_{\max}(\Sigma)$ are smaller. In turn, this means that matrices $\mathscr{A}$ and $\bm{\mathscr{A}}$ in \eqref{eq:caligA} and \eqref{eq:caligAbold} must be well conditioned and the VAR process is stable, with eigenvalues well inside the unit circle (see Section \ref{sec:stability}). Alternatively, we can also understand \eqref{eq:Nbound} as showing that large values of $\mathscr{M}$ and small values of $\mathscr L$ indicate stronger dependency in the data, thus requiring more samples for the RE conditions to hold with high probability. 

Analyzing Theorems \ref{theoremReg} and \ref{theoremRE} we can interpret the established results as follows. As the size and dimensionality $N$, $p$ and $d$ of the problem increase, we emphasize the scale of the results and use the order notations to denote the constants. Select a number of samples  at least $N \geq \mathcal{O}(w^2(\Theta))$ and let the regularization parameter satisfy $\lambda_N \geq \mathcal{O}\left(\frac{w(\Omega_R)}{\sqrt{N}} + \frac{w^2(\Omega_R)}{N^2} \right)$. With high probability then the restricted eigenvalue condition $\frac{||Z\Delta||_2}{||\Delta||_2} \geq \sqrt{\kappa N}$ for $\Delta \in \textnormal{cone}(\Omega_E)$ holds, so that $\kappa = \mathcal{O}(1)$ is a positive constant. Moreover, the norm of the estimation error in optimization problem \eqref{eq:VARestimator} is bounded by $\|\Delta\|_2 \leq \mathcal{O}\left(\frac{w(\Omega_R)}{\sqrt{N}} + \frac{w^2(\Omega_R)}{N^2}\right)\Psi(\textnormal{cone}(\Omega_{E_j}))$. Note that the norm compatibility constant $\Psi(\textnormal{cone}(\Omega_{E_j}))$ is assumed to be the same for all $j = 1, \ldots, p$, which follows from our assumption in \eqref{eq:decompNorm}.

Consider now Theorem \ref{theoremReg} and the bound on the regularization parameter $\lambda_N \geq \mathcal{O}\left(\frac{w(\Omega_R)}{\sqrt{N}} + \frac{w^2(\Omega_R)}{N^2} \right)$. As the dimensionality of the problem $p$ and $d$ grows and the number of samples $N$ increases, the first term $\frac{w(\Omega_R)}{\sqrt{N}}$ will dominate the second one $\frac{w^2(\Omega_R)}{N^2}$. This can be seen by computing $N$ for which the two terms become equal $\frac{w(\Omega_R)}{\sqrt{N}} = \frac{w^2(\Omega_R)}{N^2}$, which happens at $N = w^{\frac{2}{3}}(\Omega_R) < w(\Omega_R)$. Therefore, we can rewrite our results as follows: once the restricted eigenvalue condition holds and $\lambda_N \geq \mathcal{O}\left(\frac{w(\Omega_R)}{\sqrt{N}}\right)$, the error norm is upper-bounded by $\|\Delta\|_2 \leq \mathcal{O}\left(\frac{w(\Omega_R)}{\sqrt{N}}\right)\Psi(\textnormal{cone}(\Omega_{E_j}))$.
 
\subsection{Special Cases}
\label{sec:specialCases}

While the presented results are valid for any norm $R(\cdot)$, separable along the rows of $A_k$, it is instructive to specialize our analysis to a few popular regularization choices, such as $L_1$ and Group Lasso, Sparse Group Lasso and OWL norms.  

\subsubsection{Lasso}
\label{sec:lasso}
To establish results for $L_1$ norm, we assume that the parameter $\bm{\beta}^*$ is $s$-sparse, which in our case is meant to represent the largest number of non-zero elements in any $\beta_i$, $i=1,\ldots, p$, i.e., the combined $i$-th rows of each $A_k$, $k=1,\ldots, d$. Since $L_1$ is decomposable, it can be shown that $\Psi(\textnormal{cone}(\Omega_{E_j})) \leq 4\sqrt{s}$. Next, since $\Omega_R = \{u \in \mathbb{R}^{dp}| R(u) \leq 1\}$, then using Lemma $3$ in \cite{banerjee14} and Gaussian width results in \cite{crpw12}, we can establish that $w(\Omega_{R}) \leq \mathcal{O}(\sqrt{\log(dp)})$. Therefore, based on Theorem~4.3 and the discussion at the end of Section \ref{sec:discussion}, the bound on the regularization parameter takes the form $\lambda_N \geq \mathcal{O}\left(\sqrt{\log(dp)/N}\right)$. Hence, the estimation error is bounded by $\|\Delta\|_2 \leq \mathcal{O}\left(\sqrt{s\log(dp)/N}\right)$ as long as $N > \mathcal{O}(\log(dp))$.

\subsubsection{Group Lasso}
\label{sec:glasso}

To establish results for Group norm, we assume that for each $i=1,\ldots, p$, the vector $\beta_i \in \mathbb{R}^{dp}$ can be partitioned into a set of $K$ disjoint groups, $G = \{G_1, \ldots, G_K\}$, with the size of the largest group $m = \underset{k}{\max}|G_k|$.  Group Lasso norm is defined as $\|\bm{\beta}\|_{\text{GL}} = \sum_{k=1}^K\|\beta_{G_k}\|_2$.
We assume that the parameter $\bm{\beta}^*$ is $s_G$-group-sparse, which means that the largest  number of non-zero groups in any  $\beta_i$, $i=1,\ldots, p$ is $s_G$. Since Group norm is decomposable, as was established in \cite{negahbann12}, it can be shown that $\Psi(\textnormal{cone}(\Omega_{E_j})) \leq 4\sqrt{s_G}$. Similarly as in the Lasso case, using Lemma $3$ in \cite{banerjee14}, we get $w(\Omega_{R_{\text{GL}}}) \leq \mathcal{O}(\sqrt{m + \log(K)})$. The bound on the $\lambda_N$ takes the form $\lambda_N \geq \mathcal{O}\left(\sqrt{(m+\log(K))/N}\right)$. Combining these derivations, we obtain the bound $\|\Delta\|_2 \leq \mathcal{O}\left(\sqrt{s_G(m + \log(K))/N}\right)$ for $N > \mathcal{O}(m+\log(K))$.

\subsubsection{Sparse Group Lasso}
Similarly as in Section \ref{sec:glasso}, we assume that we have $K$ disjoint groups of size at most $m$. The Sparse Group Lasso norm enforces sparsity not only across but also within the groups and is defined as $\|\bm{\beta}\|_{\text{SGL}} = \alpha\|\bm{\beta}\|_1 + (1-\alpha)\sum_{k=1}^K\|\beta_{G_k}\|_2$, where $\alpha\in [0,1]$ is a parameter which regulates a convex combination of Lasso and Group Lasso penalties. Note that since $\|\bm{\beta}\|_2 \leq \|\bm{\beta}\|_1$, it follows that $\|\bm{\beta}\|_{\text{GL}} \leq \|\bm{\beta}\|_{\text{SGL}}$. As a result, for $\bm{\beta} \in \Omega_{R_{\text{SGL}}} \Rightarrow \bm{\beta} \in \Omega_{R_{\text{GL}}}$, so that $\Omega_{R_{\text{SGL}}} \subseteq \Omega_{R_{\text{GL}}}$ and thus $w(\Omega_{R_{\text{SGL}}}) \leq w(\Omega_{R_{\text{GL}}}) \leq  \mathcal{O}(\sqrt{m + \log(K)})$, according to Section \ref{sec:glasso}. Assuming $\bm{\beta}^*$ is $s$-sparse and $s_G$-group-sparse and noting that the norm is decomposable, we get $\Psi(\textnormal{cone}(\Omega_{E_j})) \leq 4(\alpha\sqrt{s}+(1-\alpha)\sqrt{s_G}))$. Consequently, the error bound is $\|\Delta\|_2 \leq \mathcal{O}\left(\sqrt{(\alpha s +(1-\alpha) s_G)(m + \log(K))/N}\right)$.

\subsubsection{OWL norm}
Ordered weighted $L_1$ norm is a recently introduced regularizer and is defined as $\|\bm{\beta}\|_{\text{owl}} = \sum_{i=1}^{dp} c_i|\beta|_{(i)}$, where $c_1 \geq \ldots \geq c_{dp} \geq 0$ is a predefined non-increasing sequence of weights and $|\beta|_{(1)} \geq \ldots \geq |\beta|_{(dp)} $ is the sequence of absolute values of $\bm{\beta}$, ranked in decreasing order. In \cite{chba15} it was shown that $w(\Omega_{R}) \leq \mathcal{O}(\sqrt{\log(dp)}/\bar{c})$, where $\bar{c}$ is the average of $c_1, \ldots, c_{dp}$ and the norm compatibility constant is  $\Psi(\textnormal{cone}(\Omega_{E_j})) \leq 2c_1^2\sqrt{s}/\bar{c}$. Therefore,  based on Theorem~4.3, we get $\lambda_N \geq \mathcal{O}\left(\sqrt{\log(dp)/(\bar{c}N)}\right)$ and the estimation error is bounded by $\|\Delta\|_2 \leq \mathcal{O}\left(\frac{2c_1}{\bar{c}}\sqrt{s\log(dp)/(\bar{c}N)}\right)$.

We note that the bound obtained for Lasso and Group Lasso is similar to the bound obtained in \cite{song11, bami15, koca15}. Moreover, this result is also similar to the works, which dealt with independent observations, e.g., \cite{bickel09, negahbann12}, with the difference being the constants, reflecting correlation between the samples, as we discussed in Section \ref{sec:discussion}.  The explicit bound for  Sparse Group Lasso and OWL is a \emph{novel} aspect of our work for the non-asymptotic recovery guarantees for the VAR estimation problem with norm regularization, being just a simple consequence from our more general framework. 

\vspace{-4pt}
\subsection{Proof Sketch}
\label{sec:proofScetch}

In this Section we outline the steps of the proof for Theorem \ref{theoremReg} and \ref{theoremRE}, all the details can be found in Appendix \ref{sec:lambdaBound_append} and \ref{sec:re_append}.

\subsubsection{Bound on Regularization Parameter}
\label{sec:proofLambda}
Recall that our objective is to establish for $\alpha>0$ a probabilistic statement that $\lambda_N \geq \alpha \geq R^*[\frac{1}{N}Z^T\bm{\epsilon}] = \underset{R(U) \leq 1}{\text{sup}} \left<\frac{1}{N}Z^T\bm{\epsilon},  U\right>$, where $U = [u_1^T,\ldots,u_p^T]^T \in \mathbb{R}^{dp^2}$ for $u_j \in \mathbb{R}^{dp}$ and $\bm{\epsilon} = \text{vec}(E)$ for $E$ in \eqref{eq:VARmatrix}. We denote $E_{:,j} \in \mathbb{R}^{N}$ as a column of noise matrix $E$ and note that since $Z = I_{p\times p}\otimes X$, then using the row-wise separability assumption in \eqref{eq:decompNorm} we can split the overall probability statement into $p$ parts, which are easier to work with. Thus, our objective would be to establish
\begin{align}
\label{eq:objectiveInequlity}
\mathbb{P}\bigg[ \underset{R(u_j) \leq r_j}{\text{sup}}\frac{1}{N}\left<X^TE_{:,j},u_j\right> \leq \alpha_j \bigg] \geq \pi_j,
\end{align}
for $j=1,\ldots,p$, where $\sum_{j=1}^p\alpha_j = \alpha$ and $\sum_{j=1}^pr_j=1$.

The overall strategy is to first show that the random variable $\frac{1}{N}\left<X^TE_{:,j},u_j\right>$ has sub-exponential tails. Based on the generic chaining argument, we then use Theorem~1.2.7 from \cite{talagrand06} and bound the expectation $\mathbb{E}\Bigg[\underset{R(u_j)\leq r_j}{\sup}\frac{1}{N}\left<X^TE_{:,j},u_j\right>\Bigg]$. Finally, using Theorem~1.2.9 in \cite{talagrand06} we establish the high probability bound on concentration of $\underset{R(u_j)\leq r_j}{\sup}\frac{1}{N}\left<X^TE_{:,j},u_j\right>$ around its mean, i.e., derive the bound in \eqref{eq:objectiveInequlity}.

We note that the main difficulty of working with the term $\left<X^TE_{:,j},u_j\right>$ is the complicated dependency between $X$ and $E_{:,j}$, which is due to the VAR generation process in \eqref{eq:VARmatrix}. However, if we write $\left<X^TE_{:,j},u_j\right> = \sum_{i=1}^NE_{i,j}, (X_{i,:}u_j) = \sum_{i=1}^Nm_i$, where $m_i = E_{i,j}(X_{i,:}u_j)$ and we can interpret this as a summation over martingale difference sequence \cite{lutkepohl07}. This can be easily proven by showing $\mathbb{E}(m_i|m_1,\ldots,m_{i-1})=0$. The latter is true since  in $m_i = E_{i,j}(X_{i,:}u_j)$ the terms $E_{i,j}$ and $X_{i,:}u_j$ are independent since $\epsilon_{d+i}$ is independent from $x_{d-k+i}$ for $0\leq i\leq T-d$ and $1\leq k \leq d$ (see \eqref{eq:vard}).

To show that $\sum_{i=1}^NE_{i,j}, (X_{:,i}u_j)$ has sub-exponential tails, recall that since $\epsilon_t$ in \eqref{eq:vard} is Gaussian, $E_{i,j}$ and $X_{i,:}u_j$ are independent Gaussian random variables, whose product has sub-exponential tails. Moreover, the sum over sub-exponential martingale difference sequence can be shown to be itself sub-exponential using \cite{sham11}, based on Bernstein-type inequality \cite{vershynin10}.


\subsubsection{Restricted Eigenvalue Condition}
To show $\frac{||(I_{p\times p}\otimes X)\Delta||_2}{||\Delta||_2} \geq 0$ for all $\Delta \in \textnormal{cone}(\Omega_E)$, similarly as before, we split the problem into $p$ parts by using row-wise separability assumption of the norm in \eqref{eq:decompNorm}. In particular, denote $\Delta = [\Delta_1^T, \ldots, \Delta_p^T]^T$, where $\Delta_j$ is $dp\times 1$, then we can represent the original set $\Omega_E$ as a Cartesian product of subsets $\Omega_{E_j}$, i.e., $\Omega_E = \Omega_{E_1} \times \cdots \times \Omega_{E_p}$, implying that $\textnormal{cone}(\Omega_E) = \textnormal{cone}(\Omega_{E_1})\times \cdots \times \textnormal{cone}(\Omega_{E_p})$. Therefore, our objective would be to establish
\begin{align}
\label{eq:REboundType}
\mathbb{P}\Bigg[ \underset{u_j \in \Theta_j}{\inf}||Xu_j||_2\geq \nu_j\Bigg] \geq \pi_j,
\end{align}
for each $j=1,\ldots, p$, where $\Theta = \textnormal{cone}(\Omega_{E_j})\cap S^{dp-1}$ and we defined $u_j = \frac{\Delta_j}{||\Delta_j||_2}$, since it will be easier to operate with unit-norm vectors. In the following, to reduce clutter, we drop the index $j$ from the notations.

The overall strategy is to first show that $\|Xu\|_2 - \mathbb{E}(\|Xu\|_2)$ is a sub-Gaussian random variable. Then, using generic chaining argument in \cite{talagrand06}, specifically Theorem~2.1.5, we bound $\mathbb{E}\left(\underset{u \in \Theta}{\inf}||Xu||_2\right)$. Finally, based on Lemma~2.1.3 in \cite{talagrand06} we establish the concentration inequality on $\underset{u \in \Theta}{\inf}||Xu||_2$ around its mean, i.e., derive the bound in \eqref{eq:REboundType}.




\section{Experimental Results}
\label{sec:results}
In this Section we present the experiments on simulated and real data to demonstrate the obtained theoretical results. In particular, for $L_1$ and Group $L_1$, Sparse Group $L_1$ and OWL we investigate how error norm $\|\Delta\|_2$ and regularization parameter $\lambda_N$ scale as the problem size $p$ and $N$ change. Moreover, using flight data we also compare the performance of the regularizers in real world scenario.

\begin{figure*} 
\begin{center}
\centerline{\includegraphics[width=0.97\textwidth]{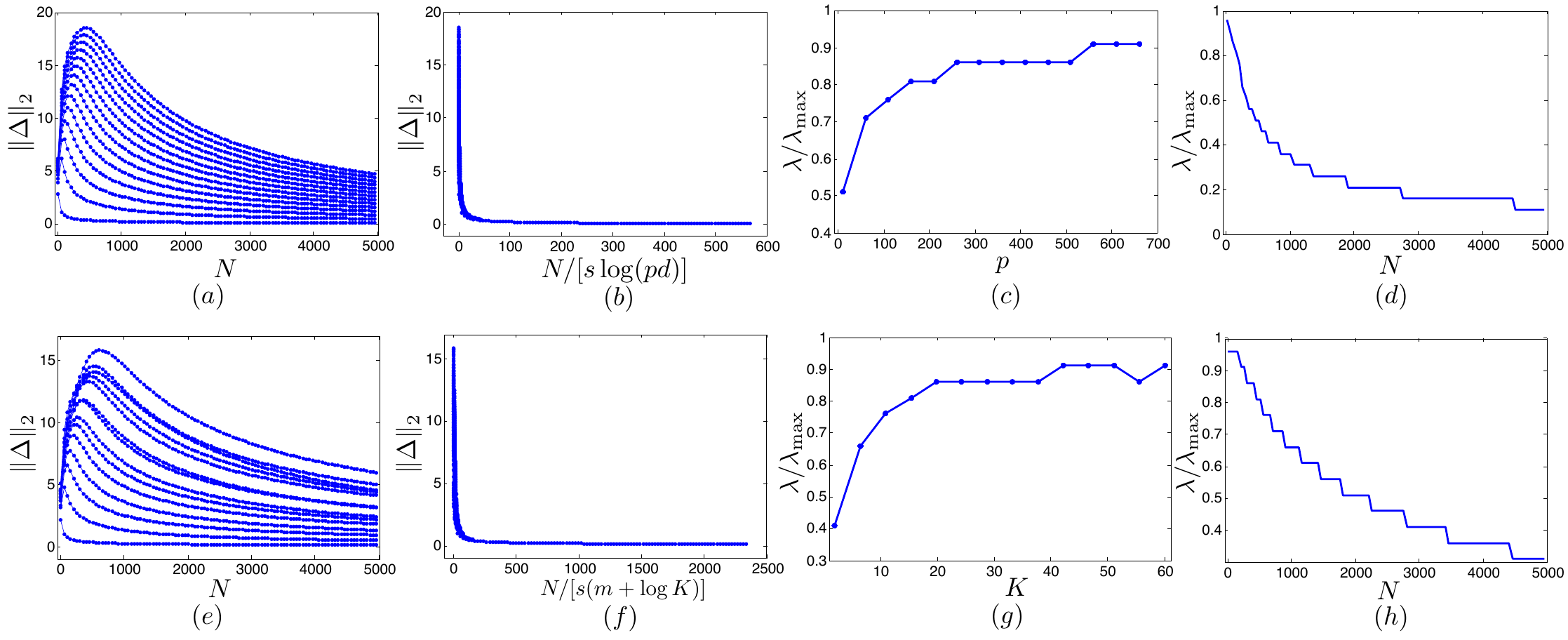}}
\vspace{-10pt}
\caption{Results for estimating parameters of a stable first order sparse VAR (top row) and group sparse VAR (bottom row). Problem dimensions: $p \in [10, 600]$, $N \in [10,5000]$, $\frac{\lambda_N}{\lambda_{max}} \in [0,1]$, $K \in [2, 60]$ and $d=1$. Figures $(a)$ and $(e)$ show dependency of errors on sample size for different $p$; in Figure $(b)$  the $N$ is scaled by $(s\log p)$ and plotted against $\|\Delta\|_2$ to show that errors scale as $(s\log p)/N$; in $(f)$ the graph is similar to $(b)$ but for group sparse VAR; in $(c)$ and $(g)$ we show dependency of $\lambda_N$ on $p$ (or number of groups $K$ in $(g)$) for fixed sample size $N$; finally, Figures $(d)$ and $(h)$ display the dependency of $\lambda_N$ on $N$ for fixed $p$.}
\label{fig:varglasso}
\end{center}
\end{figure*}

\begin{figure}
\begin{center}
\includegraphics[width=\textwidth]{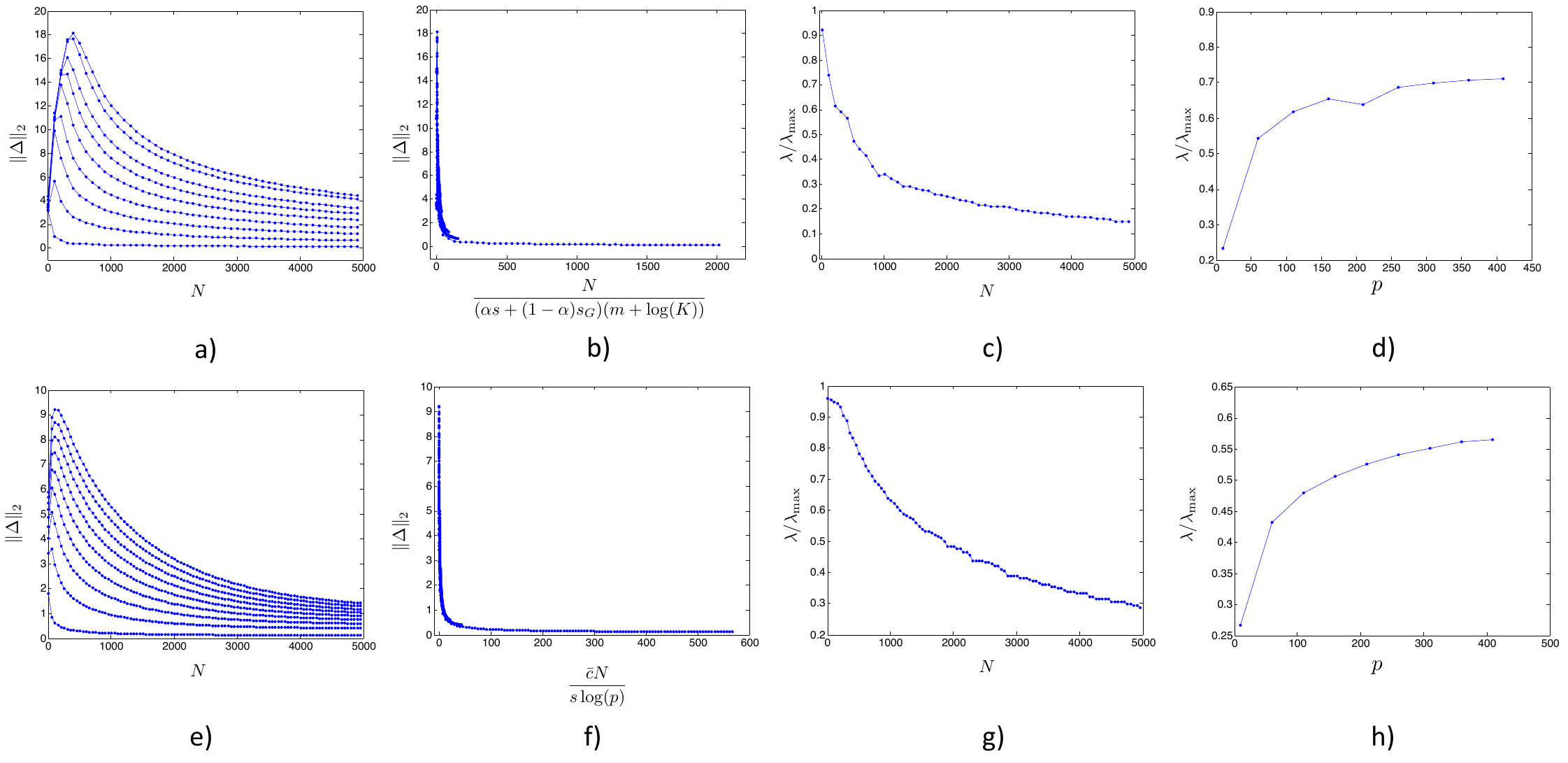}
\vspace{-25pt}
\caption{Results for estimating parameters of a stable first order Sparse Group Lasso VAR (top row) and OWL-regularized VAR (bottom row). Problem dimensions for Sparse Group Lasso : $p \in [10, 410]$, $N \in [10,5000]$, $\frac{\lambda_N}{\lambda_{max}} \in [0,1]$, $K \in [2, 60]$ and $d=1$. Problem dimensions for OWL: $p \in [10, 410]$, $N \in [10,5000]$, $\frac{\lambda_N}{\lambda_{max}} \in [0,1]$, $s \in [4, 260]$ and $d=1$. All results are shown after averaging across $50$ runs.}
\label{fig:sgl_lasso}
\end{center}
\end{figure}

\subsection{Synthetic Data}

Using synthetically generated datasets we evaluate the obtained theoretical bounds for estimation VAR under Lasso, Sparse Group Lasso, OWL and Group Lasso regularizations.
\subsubsection{Lasso}

To evaluate the estimation problem with $L_1$ norm, we simulated a first-order VAR process for different values of $p \in [10, 600]$, $s \in [4, 260]$, and $N \in [10, 5000]$. Regularization parameter was varied in the range  $\lambda_N \in (0, \lambda_{\max})$, where $\lambda_{\max}$ is the largest parameter, for which estimation problem \eqref{eq:VARestimator} produces a zero solution. All the results are shown after averaging across $50$ runs.

The results for Lasso are shown in the top row of Figure \ref{fig:varglasso}. In particular, in Figure \ref{fig:varglasso}.$a$ we show $\|\Delta\|_2$ for different $p$ and $N$ for fixed $\lambda_N$. When $N$ is small, the estimation error is large and the results cannot be trusted. However, once $N \geq \mathcal{O}(w^2(\Theta))$, the RE condition in Lemma \ref{lemmaRE} is satisfied and we see a fast decrease of errors for all $p$'s. In Figure \ref{fig:varglasso}.$b$ we plot $\|\Delta\|_2$ against rescaled sample size $\frac{N}{s\log(pd)}$. The errors are now closely aligned, confirming results of Section \ref{sec:lasso}, i.e, $\|\Delta\|_2 \leq \mathcal{O}\left(\sqrt{(s\log(pd))/N}\right)$.

Finally, in Figures \ref{fig:varglasso}.$c$ and \ref{fig:varglasso}.$d$ we show the dependence of optimal $\lambda_N$ (for fixed $N$ and $p$, we picked $\lambda_N$ achieving the smallest estimation error) on $N$ and $p$. It can be seen that as $p$ increases, $\lambda_N$ grows (for fixed $N$) at the rate similar to $\sqrt{\log p}$. On the other hand, as $N$ increases, the selected $\lambda_N$ decreases (for fixed $p$) at the rate similar to $1/\sqrt{N}$ .

\subsubsection{Sparse Group Lasso}

To evaluate the estimation problem with Sparse Group Lasso norm, we constructed first-order VAR process for the following set of problem sizes $p \in [10, 400]$, $s \in [10, 200]$, $s_G \in [2,20]$ and $N \in [10, 5000]$. The parameter $\alpha$ was set to $0.5$. Results are shown in Figure \ref{fig:sgl_lasso}, top row. Similarly as in main paper, we can see that the errors are scaled by $\frac{N}{(\alpha s + (1-\alpha)s_G)(m+\log(K))}$. Moreover, the $\lambda_N$ parameter is decreasing when number of samples $N$ increases. On the other hand, as the problem dimension $p$ increases, the selected $\lambda_N$ grows at the rate similar to $\sqrt{\log p}$.

\subsubsection{OWL}

To test the VAR estimation problem under OWL norm we constructed a first-order VAR process with $p \in [10, 410]$, $s \in [4, 260]$ and $N \in [10, 5000]$. The vector of weights $c$ was set to be a monotonically decreasing sequence of numbers in the range $[1,0)$. Figure \ref{fig:sgl_lasso}, bottom row, shows the results. It can be seen from Figure \ref{fig:sgl_lasso}-f that when the errors are plotted against $\frac{\bar{c}N}{s\log(p)}$, they become tightly aligned, confirming the bounds established in Section 3.3.4 in the main paper for the error norm. As shown in Figure \ref{fig:sgl_lasso}-g,h the selected regularization parameter $\lambda_N$ grows with the problem dimension $p$ and decreases with the number of samples $N$

\subsubsection{Group Lasso}
For Group Lasso the sparsity in rows of $A_1$ was generated in groups, whose number varied as $K \in [2,60]$. We set the largest number of non-zero groups in any row as $s_G \in [2,22]$. Results are shown in the bottom row of Figure \ref{fig:varglasso}, which have similar flavor as in Lasso case. The difference can be seen in Figure \ref{fig:varglasso}.$f$, where a close alignment of errors occurs when $N$ is now scaled as $\frac{N}{s_G(m+\log(K))}$. Moreover, the selected regularization parameter $\lambda$ increases with the number of groups $K$ and decreases with $N$. 

\begin{table*}[!tb]
\centering 
\begin{tabular}{|c|c|c|c|c|}
\hline
Lasso & OWL & Group Lasso &  Sparse Group Lasso & Ridge\\
\hline
 32.3(6.5) & 32.2(6.6) & 32.7(6.5) & 32.2(6.4) & 33.5(6.1)\\
\hline
 32.7(7.9) & 44.5(15.6) & 75.3(8.4) & 38.4(9.6) & 99.9(0.2)\\
\hline
\includegraphics[width=0.17\textwidth]{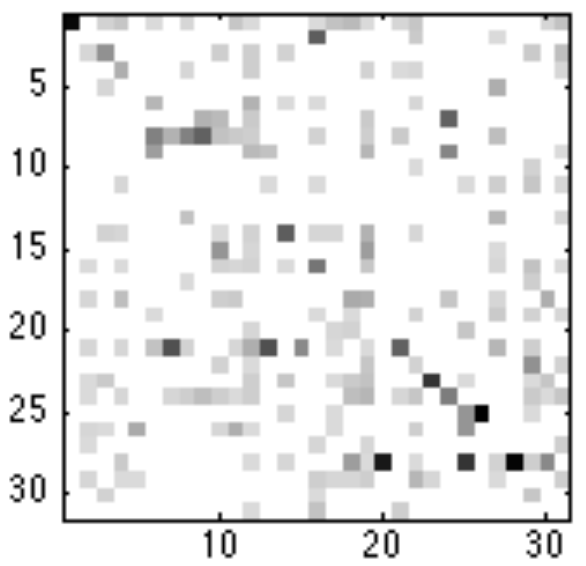} &
\includegraphics[width=0.17\textwidth]{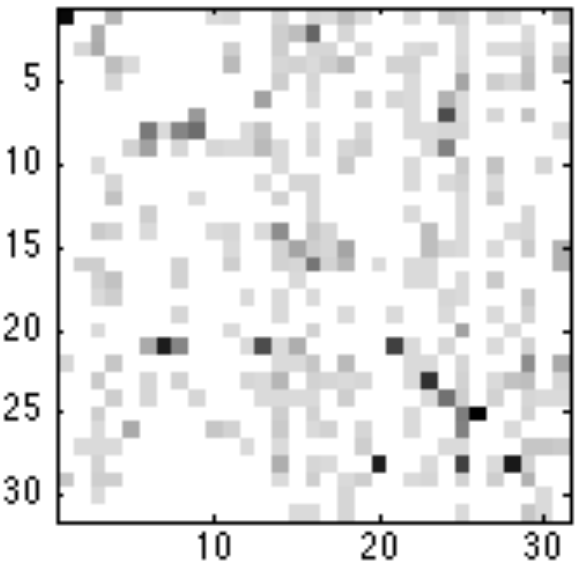} &
\includegraphics[width=0.17\textwidth]{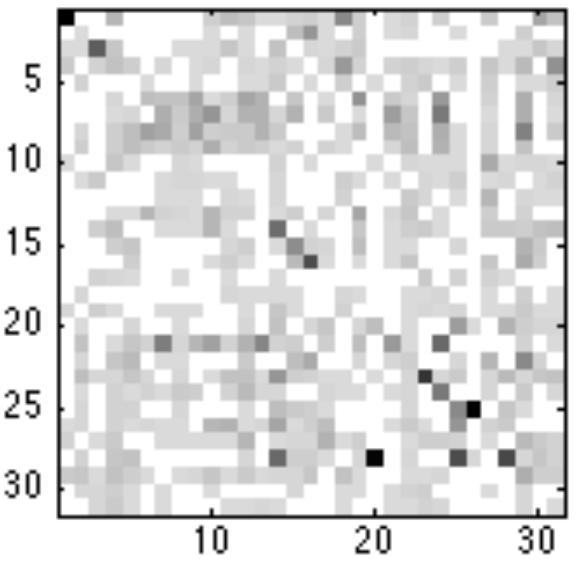} &
\includegraphics[width=0.17\textwidth]{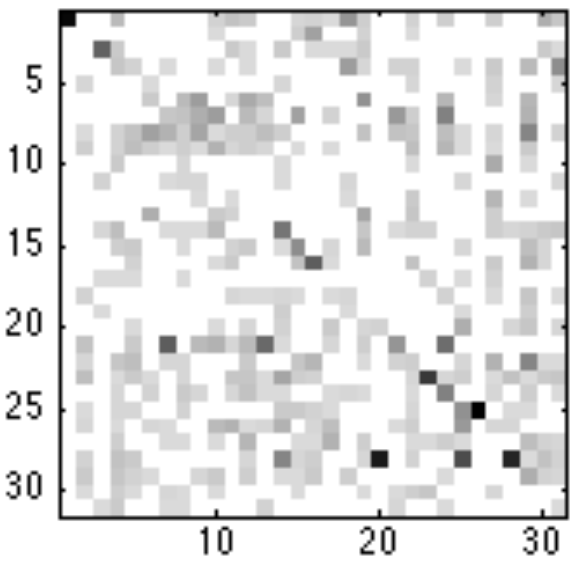} &
\includegraphics[width=0.17\textwidth]{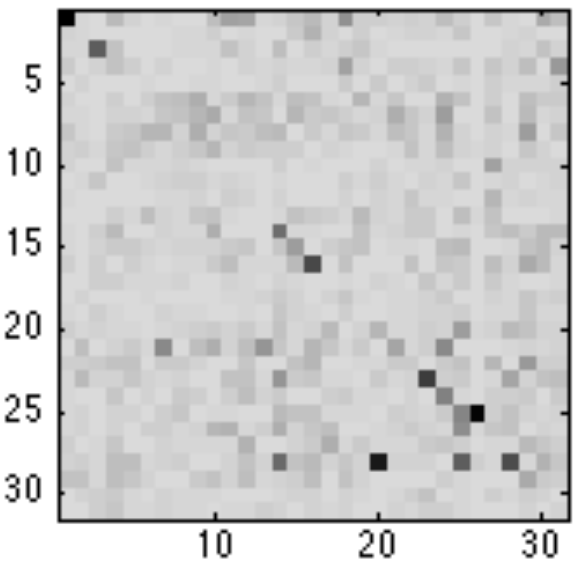} \\
\hline
\end{tabular}
\caption{Mean squared error (row 2) of the five methods used in fitting VAR model, evaluated on aviation dataset (MSE is computed using one-step-ahead prediction errors). Row 3 shows the average number of non-zeros (as a percentage of total number of elements) in the VAR matrix. The last row shows a typical sparsity pattern in $A_1$ for each method (darker dots - stronger dependencies, lighter dots - weaker dependencies). The values in parenthesis denote one standard deviation after averaging the results over 300 flights.}
\label{eq:table}
\end{table*}

\begin{table}
\centering
\begin{tabular}{|c|c|}
\hline
1 & Altitude\\
\hline
2 & Corrected angle of attack\\
\hline
3 & Brake temperature\\
\hline
4 & Computed airspeed\\
\hline
5 & Drift angle\\
\hline
6 & Engine temperature\\
\hline
7 & Low rotor speed\\
\hline
8 & High rotor speed\\
\hline
9 & Engine oil pressure\\
\hline
10 & Engine oil quantity\\
\hline
11 & Engine oil temperature\\
\hline
12 & Engine pre-cooler outlet temperature\\
\hline
13 & Fuel mass flow rate\\
\hline
14 & Lateral acceleration\\
\hline
15 & Longitudinal acceleration\\
\hline
16 & Normal acceleration\\
\hline
17 & Glide slope deviation\\
\hline
18 & Ground speed\\
\hline
19 & Localization deviation\\
\hline
20 & Magnetic heading\\
\hline
21 & Burner pressure\\
\hline
22 & Pitch angle\\
\hline
23 & Roll angle\\
\hline
24 & HPC exit temperature\\
\hline
25 & Angle magnitude\\
\hline
26 & Angle true\\
\hline
27 & Total fuel quantity\\
\hline
28 & True heading\\
\hline
29 & Vertical speed\\
\hline
30 & True airspeed\\
\hline
31 & MACH\\
\hline
\end{tabular}
\caption{ 31 features selected for structured VAR estimation on real flight data.}
\label{eq:table2}
\end{table}

\subsection{Real Data}

We have also performed evaluation tests on real data to compare the accuracy of the VAR estimation using various penalized formulations based on five norms: $L_1$, OWL, Group, Sparse Group and Ridge (square of $L_2$). Although $\|\cdot\|_2^2$ is not a norm, we included its results for reference purposes as it is frequently used in practice. In terms of data, we used the NASA flight dataset from \cite{nasadata}, consisting of over 100,000 flights, each having a record of about 250 parameters, sampled at 1 Hz. For our test, we selected 300 flights and picked 31 parameters most suitable for the prediction task (shown in Table \ref{eq:table2}) and focused on the landing part of the trajectory (duration approximately 15 minutes). For each flight we separately fitted a first-order VAR model using five approaches and performed $5$-fold cross validation to select $\lambda$, achieving smallest prediction error. For Sparse Group we set $\alpha=0.5$, while for OWL the weights $c_1,\ldots, c_p$ were set as a monotonically decreasing sequence. Table \ref{eq:table} shows the results after averaging across 300 flights.

From the table we can see that the considered problem exhibits a sparse structure since all the methods detected similar patterns in matrix $A_1$. In particular, the analysis of such patterns revealed a meaningful relationship among the flight parameters (darker dots), e.g., normal acceleration had high dependency on vertical speed and angle-of-attack, the altitude had mainly dependency with fuel quantity, vertical speed with aircraft nose pitch angle, etc.  The results also showed that the sparse regularization helps in recovering more accurate and parsimonious models as is evident by comparing performance of Ridge regression with other methods. Moreover, while all the four Lasso-based approaches performed similar to each other, their sparsity levels were different, with Lasso producing the sparsest solutions. As was also expected, Group Lasso had larger number of non-zeros since it did not enforce sparsity within the groups, as compared to the sparse version of this norm.

\section{Conclusions}
\label{sec:conc}
In this work we present a set of results for characterizing non-asymptotic estimation
error in estimating structured vector autoregressive models. The analysis holds for \emph{any} norms, separable along the rows of parameter matrices.
Our analysis is general as it is expressed in terms of Gaussian widths, a geometric measure of size of suitable sets, and includes as special cases many of the existing results focused on structured sparsity in VAR models. 



\section*{Appendix}

\appendix
\section{Stability of VAR Model}
\label{sec:stability_append}
A VAR process is stable if all the eigenvalues of $\mathbf{A}$, defined in \eqref{eq:var1}, are smaller than $1$, i.e., eigenvalues of $\mathbf{A}$ must satisfy $\text{det}(\lambda I_{dp\times dp}- \mathbf{A}) = 0$ for $\lambda \in \mathbb{C}$, $|\lambda| < 1$, $|\lambda| \neq 0$. Specifically, write
\begin{align*}
\lambda I_{dp\times dp}- \mathbf{A} &= 
\begin{bmatrix}
I\lambda & 0 & \ldots & 0 & 0\\
0 & I\lambda & \ldots & 0 & 0\\
0 & 0 & \ldots & 0 & 0\\
\vdots & \vdots & \ldots & \vdots & \vdots\\
0 & 0 & \ldots & 0 & I\lambda\\
\end{bmatrix} - 
\begin{bmatrix}
A_1 & A_2 & \ldots & A_{d-1} & A_d\\
I & 0 & \ldots & 0 & 0\\
0 & I & \ldots & 0 & 0\\
\vdots& \vdots & \ddots & \vdots & \vdots\\
0 & 0 & \ldots & I & 0
\end{bmatrix}\\
&=
\begin{bmatrix}
I\lambda - A_1 & -A_2 & \ldots & -A_{d-1} & -A_d\\
-I & I\lambda & \ldots & 0 & 0\\
0 & -I & \ldots & 0 & 0\\
\vdots& \vdots & \ddots & \vdots & \vdots\\
0 & 0 & \ldots & -I & I\lambda
\end{bmatrix}.
\end{align*}

Now multiply last ($d$-th) block-column by $\frac{1}{\lambda}$ and add to $(d-1)$-st block-column. Next, multiply the result in $(d-1)$-st block-column by $\frac{1}{\lambda}$ and add to $(d-2)$-nd block-column. Continuing in this manner, we will arrive at
\begin{align*}
Q = 
\begin{bmatrix}
\lambda I_{p\times p} - A_1 - \frac{1}{\lambda}A_2 - \ldots - \frac{1}{\lambda^{d-1}}A_d & M\\
0 & \lambda I_{p(d-1)\times p(d-1)}
\end{bmatrix},
\end{align*}
\noindent where matrix $M \in \mathbb{R}^{p\times p(d-1)}$ denotes the result of some of the column operations. Since such column operations leave the matrix determinant unchanged, we have
\begin{align*}
\text{det}(\lambda I_{dp\times dp}- \mathbf{A}) = \text{det}(Q) &= \text{det}(\lambda I_{p\times p} - A_1 - \frac{1}{\lambda}A_2 - \ldots - \frac{1}{\lambda^{d-1}}A_d)\cdot\text{det}(\lambda I_{p(d-1)\times p(d-1)})\\
& = \text{det}(I_{p\times p} - \frac{1}{\lambda}A_1 - \frac{1}{\lambda^2}A_2 - \ldots - \frac{1}{\lambda^{d}}A_d)\cdot\lambda^{pd}.
\end{align*}

Therefore, stability of VAR model in \eqref{eq:var1} requires $\text{det}(I-\sum_{k=1}^dA_k\frac{1}{\lambda^k}) = 0$ to be satisfied for $|\lambda| < 1$, $|\lambda| \neq 0$. Equivalently, $\text{det}(I-\sum_{k=1}^dA_kz^k) = 0$ must be satisfied for $z \in \mathbb{C}$, $|z| > 1$, or $\text{det}(I-\sum_{k=1}^dA_kz^k) \neq 0$ must hold for $|z| \leq 1$.

\section{Properties of Data Matrix $X$}
\label{sec:covariance_append}
In this Section we provide additional details about the covariance structure of VAR matrix $X$ as was originally presented in Section \ref{sec:Xproperties}. Recall that our VAR process is defined as
\begin{align}
\label{eq:VARmodel2}
x_t = A_1x_{t-1} + \ldots + A_dx_{t-d} + \epsilon_t, \quad t=0, \pm 1, \pm 2, \ldots,
\end{align}
where noise $\epsilon_t$ follows a Gaussian distribution, i.e., $\epsilon_t \sim \mathcal{N}(0,\Sigma)$, moreover, the distribution of $x_t$ is a zero-mean Gaussian, i.e., $x_t \sim \mathcal{N}(0, \Gamma(0))$, where $\Gamma(h) = \mathbb{E}(x_tx_{t+h}^T)$. 

Now consider the noise and data matrices from the formulation \eqref{eq:VARmatrix}
\begin{align}
\label{eq:matX}
X =
\begin{bmatrix}
x_{d-1}^T & x_{d-2}^T & \hdots &  x_{0}^T \\
x_{d}^T & x_{d-1}^T & \hdots &  x_{1}^T \\
\vdots & \vdots &  \ddots & \vdots \\
x_{T-2}^T & x_{T-3}^T & \hdots & x_{T-d-1}^T\\
x_{T-1}^T & x_{T-2}^T & \hdots & x_{T-d}^T
\end{bmatrix}.
\end{align}
In this Section our objective is to establish the probability distribution of rows of $X$.

\subsection{Single row of $X$}
\label{sec:VARcovx}
The autocovariance matrix of the original VAR process of order $d$ in \eqref{eq:VARmodel2} is defined as $\Gamma(h) = \mathbb{E}[x_tx_{t+h}^T]$. Fourier transform of autocovariance matrix is called spectral density and is denoted as (for $i = \sqrt{-1}$)
\begin{align}
\rho(\omega) = \sum_{h=-\infty}^{\infty} \Gamma(h) e^{-hi\omega}, \quad \omega \in [0, 2\pi].
\end{align}
Inverse Fourier transform of the spectral density gives back the autocovariance matrix:
\begin{align}
\Gamma(h) = \frac{1}{2\pi}\int_{0}^{2\pi} \rho(\omega) e^{hi\omega}d\omega, \quad h \in 0, \pm 1, \pm 2, \ldots
\end{align}
For our VAR model in \eqref{eq:VARmodel2}, the spectral density has a closed form expression \cite{priestley81}
\begin{align}
\label{eq:VARspectral}
\rho(\omega) = \left(I-\sum_{k=1}^dA_ke^{-ki\omega}\right)^{-1}\Sigma\left[\left(I-\sum_{k=1}^dA_ke^{-ki\omega}\right)^{-1}\right]^{*} \in \mathbb{R}^{p\times p}, 
\end{align}
where $*$ is the Hermitian of a matrix.

Let $X_{i,:}$ be any row vector of matrix $X$ in \eqref{eq:matX}, then 
\begin{align}
C_{\mathsf{X}} = 
\begin{bmatrix}
\Gamma(0) & \Gamma(1) & \ldots & \Gamma(d-1)\\
\Gamma(1)^T & \Gamma(0) & \ldots & \Gamma(d-2)\\
\vdots & \vdots & \ddots & \vdots\\
\Gamma(d-1)^T & \Gamma(d-2)^T & \ldots & \Gamma(0)\\
\end{bmatrix} \in \mathbb{R}^{dp\times dp}.
\end{align}
Note that $C_{\mathsf{X}}$ is a block-Toeplitz matrix and so we can use the following property \cite{blockToeplitz11} 
\begin{align}
\underset{
\begin{subarray}{c}
\vspace{.2em}
  1 \leq j \leq p\\
  \omega \in [0, 2\pi]
\end{subarray}}{\inf} \Lambda_j[\rho(\omega)] \leq \Lambda_k[C_V] \leq
\underset{
\begin{subarray}{c}
\vspace{.2em}
  1 \leq j \leq p\\
  \omega \in [0, 2\pi]
\end{subarray}}{\sup} \Lambda_j[\rho(\omega)], \quad \text{for  } 1 \leq k \leq Kp.
\end{align}
Using \eqref{eq:VARspectral}, we can compute the lower bound. For this we use the following relationships: for any $M$, $||M||_2 = \sqrt{\Lambda_{\text{max}}(M^TM)}$, and if $M$ is symmetric, $||M||_2 = \Lambda_{\text{max}}(M)$. Similarly, for any nonsingular $M$, $||M^{-1}||_2 = \frac{1}{\sqrt{\Lambda_{\text{min}}(M^TM)}}$, and if $M$ is symmetric, $||M^{-1}||_2 = \frac{1}{\Lambda_{\text{min}}(M)}$. Since $\rho(\omega)$ is symmetric, we have
\begin{align}
\Lambda_{\max}[\rho(\omega)] &= \left|\left|\left(I-\sum_{k=1}^dA_ke^{-ki\omega}\right)^{-1}\Sigma\left[\left(I-\sum_{k=1}^dA_ke^{-ki\omega}\right)^{-1}\right]^{*}\right|\right|_2\nonumber\\
&\leq \left|\left|\left(I-\sum_{k=1}^dA_ke^{-ki\omega}\right)^{-1}\right|\right|_2^2 ||\Sigma||_2 \nonumber\\ 
& \leq \frac{\Lambda_{\text{max}}(\Sigma)}{\Lambda_{\text{min}}\left[\left(I-\sum_{k=1}^dA_k^Te^{ki\omega}\right)\left(I-\sum_{k=1}^dA_ke^{-ki\omega}\right)\right]}
\end{align}
and the upper bound
\begin{align}
\Lambda_{\min}[\rho(\omega)] &= \left[\left|\left|\left\{\left(I-\sum_{k=1}^dA_ke^{-ki\omega}\right)^{-1}\Sigma\left[\left(I-\sum_{k=1}^dA_ke^{-ki\omega}\right)^{-1}\right]^{*}\right\}^{-1}\right|\right|_2\right]^{-1}\nonumber\\
&\geq \left[\left|\left|I-\sum_{k=1}^dA_ke^{-ki\omega}\right|\right|_2^2 ||\Sigma^{-1}||_2\right]^{-1} \nonumber\\ 
& \geq \frac{\Lambda_{\text{min}}(\Sigma)}{\Lambda_{\text{max}}\left[\left(I-\sum_{k=1}^dA_k^Te^{ki\omega}\right)\left(I-\sum_{k=1}^dA_ke^{-ki\omega}\right)\right]}.
\end{align}
Therefore, the $C_{\mathsf{X}}$ has the following bounds on its eigenvalues
\begin{align*}
\frac{\Lambda_{\text{min}}(\Sigma)}{\Lambda_{\max}\left[\left(I-\sum_{k=1}^dA_k^Te^{ki\omega}\right)\left(I-\sum_{k=1}^dA_ke^{-ki\omega}\right)\right]}
&\leq \Lambda_k[C_{\mathsf{X}}] \leq \\
&\hspace{-20pt}\leq \frac{\Lambda_{\text{max}}(\Sigma)}{\Lambda_{\min}\left[\left(I-\sum_{k=1}^dA_k^Te^{ki\omega}\right)\left(I-\sum_{k=1}^dA_ke^{-ki\omega}\right)\right]},
\end{align*}
for   $1 \leq k \leq dp, \text{ and } \omega \in [0, 2\pi]$. 

Denoting $\Lambda_{\min}(\mathscr{A}) = \Lambda_{\min}\left[\left(I-\sum_{k=1}^dA_k^Te^{ki\omega}\right)\left(I-\sum_{k=1}^dA_ke^{-ki\omega}\right)\right]$ for $\omega \in [0, 2\pi]$ and similarly $\Lambda_{\max}(\mathscr{A}) = \Lambda_{\max}\left[\left(I-\sum_{k=1}^dA_k^Te^{ki\omega}\right)\left(I-\sum_{k=1}^dA_ke^{-ki\omega}\right)\right]$ for $\omega \in [0, 2\pi]$, we can compactly write the above as 
\begin{align}
\label{eq:autocovBounds}
\frac{\Lambda_{\min}(\Sigma)}{\Lambda_{\max}(\mathscr{A})}
\leq \Lambda_k[C_{\mathsf{X}}] \leq
\frac{\Lambda_{\max}(\Sigma)}{\Lambda_{\min}(\mathscr{A})},
\end{align}
for   $1 \leq k \leq dp$. From the above we extract the lower bound and denote it as
\begin{align}
\Lambda_k[C_{\mathsf{X}}] \geq \frac{\Lambda_{\min}(\Sigma)}{\Lambda_{\max}(\mathscr{A})} = \mathscr{L}.
\end{align}

\subsection{All the rows of $X$}
\label{sec:VARcovX}
Consider a model obtained from the rows of matrix $X$ (see \eqref{eq:matX}), i.e., 
\begin{align*}
\begin{bmatrix}
x_{d-i+1}\\
x_{d-i}\\
\vdots\\
x_{i}
\end{bmatrix} =
\begin{bmatrix}
A_1 & A_2 & \ldots & A_{d-1} & A_d\\
I & 0 & \ldots & 0 & 0\\
0 & I & \ldots & 0 & 0\\
\vdots& \vdots & \ddots & \vdots & \vdots\\
0 & 0 & \ldots & I & 0
\end{bmatrix} 
\begin{bmatrix}
x_{d-i}\\
x_{d-i-1}\\
\vdots\\
x_{i-1}
\end{bmatrix} + 
\begin{bmatrix}
\epsilon_{d-i+1}\\
0\\
\vdots\\
0
\end{bmatrix}.
\end{align*}
Written in a compact form, the above expression takes the form
\begin{align*}
X_{j,:} = \mathbf{A}X_{j-1,:} + \mathcal{E}_j, \quad \text{ for }  j=1,\ldots,N,
\end{align*}
which can be thought to be the transformations of the form
\begin{align*}
X_{1,:} = 
\begin{bmatrix}
x_{d-1}\\
x_{d-2}\\
\vdots\\
x_{0}
\end{bmatrix}
\quad\rightarrow\quad
X_{2,:} = 
\begin{bmatrix}
x_{d}\\
x_{d-1}\\
\vdots\\
x_{1}
\end{bmatrix}
\quad
\rightarrow
\cdots \rightarrow
\quad
X_{N,:} = 
\begin{bmatrix}
x_{N+d-2}\\
x_{N+d-3}\\
\vdots\\
x_{N-1}
\end{bmatrix}.
\end{align*}
Let
\begin{align}
\label{eq:U}
\mathcal{U} = 
\begin{bmatrix}
X_{1,:}\\
\vdots\\
X_{N,:}
\end{bmatrix} \in \mathbb{R}^{Ndp},
\end{align}
be a vector composed from the output of the above VAR model during $N$ steps. Then $C_\mathcal{U} \in \mathbb{R}^{Ndp\times Ndp}$ is the covariance matrix of vector $\mathcal{U}$
\begin{align}
\label{eq:covU2}
C_\mathcal{U} = 
\mathbb{E}(\mathcal{U}\mathcal{U}^T) &= 
\mathbb{E}\begin{bmatrix}
X_{1,:}\\
\vdots\\
X_{N,:}
\end{bmatrix}
\begin{bmatrix}
X_{1,:}^T \ldots X_{N,:}^T
\end{bmatrix}
=\begin{bmatrix}
\mathbb{E}[X_{1,:}X_{1,:}^T] & \mathbb{E}[X_{1,:}X_{2,:}^T] & \ldots & \mathbb{E}[X_{1,:}X_{N,:}^T]\\
\mathbb{E}[X_{2,:}X_{1,:}^T] & \mathbb{E}[X_{2,:}X_{2,:}^T] & \ldots & \mathbb{E}[X_{2,:}X_{N,:}^T]\\
\vdots & \vdots & \ddots & \vdots\\
\mathbb{E}[X_{N,:}X_{1,:}^T] & \mathbb{E}[X_{N,:}X_{2,:}^T] & \ldots & \mathbb{E}[X_{N,:}X_{N,:}^T]
\end{bmatrix}.
\end{align}
To establish the bounds on the eigenvalues of $C_\mathcal{U}$, we denote the spectral density of the corresponding VAR process as
\begin{align*}
\rho_{\mathsf{X}}(\omega) = \sum_{h=-\infty}^{\infty} \Gamma_{\mathsf{X}}(h) e^{-hi\omega}, \quad \omega \in [0, 2\pi],
\end{align*}
where  $\Gamma_{\mathsf{X}}(h) = \mathbb{E}[X_{j,:}X_{j+h,:}^T]$. Since $C_\mathcal{U}$ is a block-Toeplitz matrix, we can employ the same relationship as we used in Section \ref{sec:VARcovx} 
\begin{align}
\label{eq:covUbounds}
\underset{
\begin{subarray}{c}
\vspace{.2em}
  1 \leq l \leq dp\\
  \omega \in [0, 2\pi]
\end{subarray}}{\inf} \Lambda_l[\rho_{\mathsf{X}}(\omega)] \leq \Lambda_k[C_\mathcal{U}] \leq
\underset{
\begin{subarray}{c}
\vspace{.2em}
  1 \leq l \leq dp\\
  \omega \in [0, 2\pi]
\end{subarray}}{\sup} \Lambda_l[\rho_{\mathsf{X}}(\omega)], \quad \text{for  } 1 \leq k \leq Ndp.
\end{align}
In the following we establish the closed form expression of spectral density $\rho_{\mathsf{X}}$. For this we write
\begin{align}
\label{eq:gamma}
\rho_{\mathsf{X}}(\omega) &= \sum_{h=-\infty}^{\infty} \Gamma_{\mathsf{X}}(h) e^{-hi\omega}\nonumber\\
&=\sum_{h=-\infty}^{\infty} \mathbb{E}[X_{j,:}X_{j+h,:}^T] e^{-hi\omega} \quad \text{ for any } j\nonumber\\
&=\sum_{h=-\infty}^{\infty} \mathbb{E}\Bigg[\sum_{k=0}^{\infty}\mathbf{A}^{k}E_{j-k,:} \bigg(\sum_{s=0}^{\infty}\mathbf{A}^{s}E_{j+h-s,:}\bigg)^T\Bigg] e^{-hi\omega}\nonumber\\
&=\sum_{h=-\infty}^{\infty} \mathbb{E}\Bigg[\sum_{k=0}^{\infty}\mathbf{A}^{k}E_{j-k,:} \bigg(\sum_{s=0}^{\infty}\mathbf{A}^{s-h}E_{j-s,:}\bigg)^T\Bigg] e^{-hi\omega}\nonumber\\
&=\sum_{h=-\infty}^{\infty} \sum_{k=0}^{\infty}\mathbf{A}^{k}\Sigma_E \Big(\mathbf{A}^{k-h}\Big)^T e^{-hi\omega}\nonumber\\
&=\sum_{h=-\infty}^{\infty} \sum_{k=0}^{\infty}\mathbf{A}^{k}\Sigma_E \Big(\mathbf{A}^{k-h}\Big)^T e^{-hi\omega + ki\omega - ki\omega}\nonumber\\
&=\sum_{h=-\infty}^{\infty} \sum_{k=0}^{\infty}\mathbf{A}^{k}e^{-ki\omega}\Sigma_E \Big(\mathbf{A}^{k-h}e^{-(k-h)i\omega}\Big)^* \nonumber\\
&= \sum_{k=0}^{\infty}\mathbf{A}^{k}e^{-ki\omega}\Sigma_E \sum_{r=0}^{\infty}\Big(\mathbf{A}^{r}e^{-ri\omega}\Big)^*\nonumber\\
&= \Big(I - \mathbf{A}e^{-i\omega}\Big)^{-1}\Sigma_E \Bigg[\Big(I - \mathbf{A}e^{-i\omega}\Big)^{-1}\Bigg]^{*},
\end{align}
where we have used the fact that $\sum_{k=0}^{\infty}\mathbf{A}^{k}e^{-ki\omega} = \Big(I - \mathbf{A}e^{-i\omega}\Big)^{-1}$.

Now, using \eqref{eq:covUbounds}, \eqref{eq:gamma}, the results from Section \ref{sec:VARcovx} and the fact that the covariance matrix $\Sigma_{\mathcal{E}}$ has the form
\begin{align*}
\Sigma_{\mathcal{E}} = 
\begin{bmatrix}
\Sigma & 0 &  \ldots & 0\\
0 & 0 &  \ldots & 0\\
\vdots & \vdots & \ddots & \vdots\\
0 & 0 &  \ldots & 0\\
\end{bmatrix},
\end{align*}
we can establish the following bounds
\begin{align*}
\frac{\Lambda_{\text{min}}(\Sigma_{\mathcal{E}})}{\Lambda_{\text{max}}\left[\left(I - \mathbf{A}^Te^{i\omega}\right)\left(I - \mathbf{A}e^{-i\omega}\right)\right]} 
\leq \Lambda_k[C_\mathcal{U}] \leq
\frac{\Lambda_{\text{max}}(\Sigma_{\mathcal{E}})}{\Lambda_{\text{min}}\left[\left(I - \mathbf{A}^Te^{i\omega}\right)\left(I - \mathbf{A}e^{-i\omega}\right)\right]}.
\end{align*}
Since $\Lambda_{\max}(\Sigma_{\mathcal{E}}) = \Lambda_{\max}(\Sigma)$, the upper bound becomes
\begin{align*}
\Lambda_{\max}[C_\mathcal{U}] \leq
\frac{\Lambda_{\text{max}}(\Sigma)}{\Lambda_{\text{min}}\left[\left(I - \mathbf{A}^Te^{i\omega}\right)\left(I - \mathbf{A}e^{-i\omega}\right)\right]},
\end{align*}
for $\omega \in [0, 2\pi]$. Denoting $\Lambda_{\min}(\bm{\mathscr{A}}) = \Lambda_{\text{min}}\left[\left(I - \mathbf{A}^Te^{i\omega}\right)\left(I - \mathbf{A}e^{-i\omega}\right)\right]$ for $\omega \in [0, 2\pi]$, we can compactly write the above as 
\begin{align}
\label{eq:autocovBounds2}
\Lambda_{\max}[C_\mathcal{U}] \leq
\frac{\Lambda_{\text{max}}(\Sigma)}{\Lambda_{\min}(\bm{\mathscr{A}})}.
\end{align}

\subsection{Linear combination of rows of $X$}
\label{sec:rowDistrr}

Consider a vector $q = Xa \in \mathbb{R}^{N}$ for any $a \in \mathbb{R}^{dp}$.  Since each element $X_{i,:}^Ta \sim \mathcal{N}(0, a^TC_{\mathsf{X}}a)$, it follows that $q \sim \mathcal{N}(0, Q_a)$ with a covariance matrix $Q_a \in \mathbb{R}^{N\times N}$, which is defined as
\begin{align*}
Q_a = \mathbb{E}(qq^T) &= 
\mathbb{E}
\begin{bmatrix}
X_{1,:}^Ta\\
\vdots\\
X_{N,:}^Ta
\end{bmatrix}
\begin{bmatrix}
a^TX_{1,:} \ldots a^TX_{N,:}
\end{bmatrix}\\
&=\begin{bmatrix}
a^T\mathbb{E}[X_{1,:}X_{1,:}^T]a & a^T\mathbb{E}[X_{1,:}X_{2,:}^T]a & \ldots & a^T\mathbb{E}[X_{1,:}X_{N,:}^T]a\\
a^T\mathbb{E}[X_{2,:}X_{1,:}^T]a & a^T\mathbb{E}[X_{2,:}X_{2,:}^T]a & \ldots & a^T\mathbb{E}[X_{2,:}X_{N,:}^T]a\\
\vdots & \vdots & \ddots & \vdots\\
a^T\mathbb{E}[X_{N,:}X_{1,:}^T]a & a^T\mathbb{E}[X_{N,:}X_{2,:}^T]a & \ldots & a^T\mathbb{E}[X_{N,:}X_{N,:}^T]a
\end{bmatrix}\\
&=
\begin{bmatrix}
a^T & 0 & \ldots & 0\\
0 & a^T & \ldots & 0\\
\vdots & \vdots & \ddots & \vdots\\
0 & 0 & \ldots & a^T\\
\end{bmatrix}
\begin{bmatrix}
\mathbb{E}[X_{1,:}X_{1,:}^T] & \mathbb{E}[X_{1,:}X_{2,:}^T] & \ldots & \mathbb{E}[X_{1,:}X_{N,:}^T]\\
\mathbb{E}[X_{2,:}X_{1,:}^T] & \mathbb{E}[X_{2,:}X_{2,:}^T] & \ldots & \mathbb{E}[X_{2,:}X_{N,:}^T]\\
\vdots & \vdots & \ddots & \vdots\\
\mathbb{E}[X_{N,:}X_{1,:}^T] & \mathbb{E}[X_{N,:}X_{2,:}^T] & \ldots & \mathbb{E}[X_{N,:}X_{N,:}^T]
\end{bmatrix}
\begin{bmatrix}
a & 0 & \ldots & 0\\
0 & a & \ldots & 0\\
\vdots & \vdots & \ddots & \vdots\\
0 & 0 & \ldots & a\\
\end{bmatrix}\\
& =
(I_{N\times N}\otimes a^T)
\begin{bmatrix}
\mathbb{E}[X_{1,:}X_{1,:}^T] & \mathbb{E}[X_{1,:}X_{2,:}^T] & \ldots & \mathbb{E}[X_{1,:}X_{N,:}^T]\\
\mathbb{E}[X_{2,:}X_{1,:}^T] & \mathbb{E}[X_{2,:}X_{2,:}^T] & \ldots & \mathbb{E}[X_{2,:}X_{N,:}^T]\\
\vdots & \vdots & \ddots & \vdots\\
\mathbb{E}[X_{N,:}X_{1,:}^T] & \mathbb{E}[X_{N,:}X_{2,:}^T] & \ldots & \mathbb{E}[X_{N,:}X_{N,:}^T]
\end{bmatrix}
(I_{N\times N}\otimes a).
\end{align*}
We denote the covariance matrix in the middle as 
\begin{align}
\label{eq:covU}
C_\mathcal{U} = 
\mathbb{E}(\mathcal{U}\mathcal{U}^T) &= 
\mathbb{E}\begin{bmatrix}
X_{1,:}\\
\vdots\\
X_{N,:}
\end{bmatrix}
\begin{bmatrix}
X_{1,:}^T \ldots X_{N,:}^T
\end{bmatrix}
=\begin{bmatrix}
\mathbb{E}[X_{1,:}X_{1,:}^T] & \mathbb{E}[X_{1,:}X_{2,:}^T] & \ldots & \mathbb{E}[X_{1,:}X_{N,:}^T]\\
\mathbb{E}[X_{2,:}X_{1,:}^T] & \mathbb{E}[X_{2,:}X_{2,:}^T] & \ldots & \mathbb{E}[X_{2,:}X_{N,:}^T]\\
\vdots & \vdots & \ddots & \vdots\\
\mathbb{E}[X_{N,:}X_{1,:}^T] & \mathbb{E}[X_{N,:}X_{2,:}^T] & \ldots & \mathbb{E}[X_{N,:}X_{N,:}^T]
\end{bmatrix}.
\end{align}
Thus, we established that $q \sim \mathcal{N}(0, Q_a)$, where $Q_a = (I\otimes a^T)C_\mathcal{U}(I\otimes a)$.

In what follows, we compute $\text{trace}(Q_a)$ and $||Q_a||_2$ for the covariance matrix $Q_a$. It can be seen that the trace of $Q_a$ is given by
\begin{align}
\label{eq:Qtrace}
\text{trace}(Q_a) = Na^TC_\mathsf{X}a,
\end{align}
where $C_\mathsf{X}$ is defined in \eqref{eq:covX}. Next, we compute upper bound on $||Q_a||_2$ as follows
\begin{align}
\label{eq:Qnorm}
||Q_a||_2 &= ||(I\otimes a^T)C_\mathcal{U}(I\otimes a)||_2\nonumber\\
&\leq ||I\otimes a||_2^2~ ||C_\mathcal{U}||_2\nonumber\\
& = ||a||_2^2~\Lambda_{\max}(C_\mathcal{U}),
\end{align}
where the last equality follows since $||I\otimes a||_2^2 = \Lambda_{\max}\Big((I\otimes a^T)(I\otimes a)\Big) = \Lambda_{\max}\Big(I\otimes a^Ta\Big) = ||a||_2^2$. We used a property of Kronecker product which states that for matrices with suitable dimensions, $(A\otimes B)(C\otimes D) = (AC\otimes BD)$.

To establish $\Lambda_{\max}(C_\mathcal{U})$, we use the results from Section \ref{sec:VARcovX}, expression \eqref{eq:autocovBounds2}, which enable us to conclude that the upper bound of the largest eigenvalue of matrix $C_{\mathcal{U}}$ is given by
\begin{align*}
\Lambda_{\max}(C_{\mathcal{U}}) \leq 
\frac{\Lambda_{\max}(\Sigma)}{\Lambda_{\min}(\bm{\mathscr{A}})}.
\end{align*}
Therefore, the bound on the covariance matrix $||Q_a||_2$ in \eqref{eq:Qnorm} is now given by
\begin{align}
\label{eq:Qnorm2}
||Q_a||_2 \leq ||a||_2^2~\frac{\Lambda_{\max}(\Sigma)}{\Lambda_{\min}(\bm{\mathscr{A}})} = \mathscr{M}.
\end{align}

\section{Bound on Regularization Parameter}
\label{sec:lambdaBound_append}
To establish lower bound on the regularization parameter $\lambda_N$, we derive an upper bound on $R^*[\frac{1}{N}Z^T\bm{\epsilon}] \leq \alpha$, for some $\alpha>0$, which will establish the required relationship $\lambda_N \geq \alpha \geq R^*[\frac{1}{N}Z^T\bm{\epsilon}]$. We will also utilize the notions of Gaussian width and covering net.
\begin{definition}
\label{eq:gaussWidth}
For any set $\mathcal{S}$ and for a vector of independent zero-mean unit variance Gaussian variables $g \sim \mathcal{N}(0, I)$, the Gaussian width of the set is defined as 
\begin{align}
\label{eq:width}
w(\mathcal{S}) = \mathbb{E}_g[\underset{u\in\mathcal{S}}{\textnormal{sup}\langle g, u\rangle}].
\end{align}
\end{definition}

\noindent Denote $E_{:,j} \in \mathbb{R}^{N}$ as a column of matrix $E$ and vector $U = [u_1^T,\ldots,u_p^T]^T \in \mathbb{R}^{dp^2}$, where $u_i \in \mathbb{R}^{dp}$. Note that since $Z = I_{p\times p}\otimes X$, and $\bm{\epsilon} = \text{vec}(E)$, we can observe the following
\begin{align}
\label{eq:regBoundIneq}
\underset{R(U) \leq 1}{\text{sup}} \left<\frac{1}{N}Z^T\bm{\epsilon}, U\right> &= \underset{R(U) \leq 1}{\text{sup}} \frac{1}{N}\Bigg<\left(I_{p\times p}\otimes X^T\right)\text{vec}(E), U\Bigg>\nonumber\\
&= \underset{R([u_1^T,\ldots,u_p^T]^T) \leq 1}{\text{sup}} \frac{1}{N}\Bigg(\left<X^TE_{:,1},u_1\right> +, \ldots,+ \left<X^TE_{:,p},u_p\right>\Bigg)\nonumber\\
&= \frac{1}{N}\Bigg(\underset{R([u_1^T,\ldots,u_p^T]^T) \leq 1}{\text{sup}}\left<X^TE_{:,1},u_1\right> +, \ldots, + \underset{R([u_1^T,\ldots,u_p^T]^T) \leq 1}{\text{sup}}\left<X^TE_{:,p},u_p\right>\Bigg)\nonumber\\
&= \frac{1}{N}\Bigg(\underset{R(u_1) \leq r_1}{\text{sup}}\left<X^TE_{:,1},u_1\right> +, \ldots, +\underset{R(u_p) \leq r_p}{\text{sup}}\left<X^TE_{:,p},u_p\right>\Bigg)\nonumber\\
&= \frac{1}{N}\sum_{j=1}^p\underset{R(u_j) \leq r_j}{\text{sup}}\left<X^TE_{:,j},u_j\right> 
\end{align}
\noindent where $\sum_{j=1}^pr_j \leq 1$ and $r_j \geq 0$. 

Our objective is to establish a high probability bound of the form 
\begin{align*}
\mathbb{P}\Bigg[ \underset{R(U) \leq 1}{\text{sup}} \left<\frac{1}{N}Z^T\bm{\epsilon}, U\right> \leq \alpha \Bigg] \geq \pi
\end{align*}
where $ 0 \leq \pi \leq 1$, i.e., upper bound should hold with at least probability $\pi$. Using \eqref{eq:regBoundIneq}  and assuming that $\alpha = \sum_{j=1}^p\alpha_j$, we can rewrite the above probabilistic statement as follows
\begin{align}
\label{eq:multiStepBound}
\mathbb{P}\Bigg[ \underset{R(U) \leq 1}{\text{sup}} \left<\frac{1}{N}Z^T\bm{\epsilon}, U\right> \leq \alpha \Bigg] &= \mathbb{P}\Bigg[ \frac{1}{N}\sum_{j=1}^p\underset{R(u_j) \leq r_j}{\text{sup}}\left<X^TE_{:,j},u_j\right> \leq \sum_{j=1}^p\alpha_j \Bigg]\nonumber\\
&\geq \mathbb{P}\Bigg[ \bigg\{\underset{R(u_1) \leq r_1}{\text{sup}}\frac{1}{N}\left<X^TE_{:,1},u_1\right> \leq \alpha_1\bigg\} \text{ and }\\
 &\hspace{50pt}\ldots \text{ and }\bigg\{\underset{R(u_p) \leq r_p}{\text{sup}}\frac{1}{N}\left<X^TE_{:,p},u_p\right> \leq \alpha_p\bigg\}\Bigg]\nonumber\\
&\geq \sum_{j=1}^p \mathbb{P}\bigg[ \underset{R(u_j) \leq r_j}{\text{sup}}\frac{1}{N}\left<X^TE_{:,j},u_j\right> \leq \alpha_j \bigg] - (p-1).
\end{align}
In the above derivations we used the observation that if the events  $\bigg\{\underset{R(u_j) \leq r_j}{\text{sup}}\frac{1}{N}\left<X^TE_{:,j},u_j\right> \leq \alpha_j\bigg\}$, for each $j$ hold, then the event $\left\{\sum_{j=1}^p\underset{R(u_j) \leq r_j}{\text{sup}}\frac{1}{N}\left<X^TE_{:,j},u_j\right> \leq \sum_{j=1}^p\alpha_j\right\}$ also holds but the reverse is not always true, implying that the probability space related to the event\\ $\left\{\sum_{j=1}^p\underset{R(u_j) \leq r_j}{\text{sup}}\frac{1}{N}\left<X^TE_{:,j},u_j\right> \leq \sum_{j=1}^p\alpha_j\right\}$ is larger.

Therefore, based on \eqref{eq:multiStepBound}, we see that we need to establish the following concentration bound
\begin{align}
\label{eq:objectiveInequlityy}
\mathbb{P}\bigg[ \underset{R(u_j) \leq r_j}{\text{sup}}\frac{1}{N}\left<X^TE_{:,j},u_j\right> \leq \alpha_j \bigg] \geq \pi_j,
\end{align}
for each $j=1,\ldots,p$.

In the following our objective would be to first establish that the random variable $\frac{1}{N}\left<X^TE_{:,j},h\right>$ has sub-exponential tails, where $h \in \mathbb{R}^{dp}$, $\|h\|_2=1$ is a unit norm vector. Based on the generic chaining argument we then use Theorem 1.2.7 in \cite{talagrand06} and bound the expectation of the supremum of the original variable $\frac{1}{N}\left<X^TE_{:,j},u_j\right>$, i.e., bound $\mathbb{E}\Bigg[\underset{R(u_j)\leq r_j}{\sup}\frac{1}{N}\left<X^TE_{:,j},u_j\right>\Bigg]$. Finally, using Theorem 1.2.9 in \cite{talagrand06} we establish the high probability bound on how $\underset{R(u_j)\leq r_j}{\sup}\frac{1}{N}\left<X^TE_{:,j},u_j\right>$ concentrates around its mean.

\subsection{Martingale difference sequence}

We start by writing
\begin{align*}
\left<X^TE_{:,j},h\right> = \left<E_{:,j}, Xh\right> = \sum_{i=1}^NE_{i,j}, (X_{:,i}h) = \sum_{i=1}^Nm_i,
\end{align*}
where $m_i = E_{i,j}(X_{i,:}h)$, $i=1,\ldots,N$. Observe that $m_i$ is a martingale difference sequence (MDS), which can be shown by establishing that $\mathbb{E}(m_i|m_1,\ldots,m_{i-1})=0$ (see \cite{lutkepohl07}). We can introduce a set $\{E_{1,:}, E_{2,:}, \ldots, E_{i-1,:}\} = \{\epsilon_d^T, \epsilon_{d+1}^T, \ldots, \epsilon_T^T\}$ and write 
\begin{align*}
\mathbb{E}\big[m_i|m_1,\ldots,m_{i-1}\big] = \mathbb{E}\big[\mathbb{E}\big[m_i|m_1,\ldots,m_{i-1},E_{1,:}, \ldots, E_{i-1,:} \big]\big],
\end{align*}
using the technique of iterated expectation. Note that the set $\{E_{1,:}, E_{2,:}, \ldots, E_{i-1,:}\}$ contains more information than the set $\{m_1,\ldots,m_{i-1}\}$ and conditioning on it has fixed all the past history of the sequence until time stamp $i$. Since $m_i = E_{i,j}(X_{i,:}h)$, the terms $E_{i,j}$ and $X_{i,:}h$ are now independent. The independence follows since every row of matrix $X$ is independent of the corresponding row of matrix $E$:
\begin{align*}
E = 
\begin{bmatrix}
\epsilon_d^T \\
\epsilon_{d+1}^T\\
\vdots\\
\epsilon_{T-1}^T\\
\epsilon_T^T 
\end{bmatrix}, 
\quad
X =
\begin{bmatrix}
x_{d-1}^T & x_{d-2}^T & \hdots &  x_{0}^T \\
x_{d}^T & x_{d-1}^T & \hdots &  x_{1}^T \\
\vdots & \vdots &  \ddots & \vdots \\
x_{T-2}^T & x_{T-3}^T & \hdots & x_{T-d-1}^T\\
x_{T-1}^T & x_{T-2}^T & \hdots & x_{T-d}^T
\end{bmatrix},
\end{align*}
which can be verified by noting that the noise vector $\epsilon_{d+i}$ is independent from $x_{d-k+i}$ since $(d+i) > (d-k+i)$ for $0\leq i\leq T-d$ and $1\leq k \leq d$. In other words, the information contained in $x_{d-k+i}$ does not contain information from the noise $\epsilon_{d+i}$ (see \eqref{eq:VARmatrix}). Moreover,
\begin{align}
\label{eq:meanMi}
\mathbb{E}\Big[m_i\Big] = \mathbb{E}\Big[E_{i,j}(X_{i,:}h)\Big] = \mathbb{E}\Big[E_{i,j}\Big]\mathbb{E}\Big[X_{i,:}h\Big] = 0,
\end{align}
due to the zero-mean noise $\mathbb{E}\big[E_{i,j}\big]=0$. Consequently, we have shown that the conditional expectation $\mathbb{E}\big[m_i|m_1,\ldots,m_{i-1},E_{1,:}, \ldots, E_{i-1,:} \big]=0$ and therefore
\begin{align*}
\mathbb{E}\big[m_i|m_1,\ldots,m_{i-1}\big] = 0,
\end{align*}
proving that $m_i = E_{i,j}(X_{i,:}h)$, $i=1,\ldots,N$ is a martingale difference sequence.

Next, to show that $\frac{1}{N}\left<X^TE_{:,j},h\right> = \frac{1}{N}\sum_{i=1}^Nm_i$ has sub-exponential tails, we first show that $m_i$ is sub-exponential random variable and then use the proof argument similar to Azuma-type \cite{azuma67} and Bernstein-type \cite{vershynin10} inequalities to establish that a sum over sub-exponential martingale difference sequence is itself sub-exponential. 

\subsection{Sub-exponential tails of $\frac{1}{N}\left<X^TE_{:,j},h\right>$}

The MDS $m_i$ is sub-exponential since it is a product of two Gaussians. Indeed, recall that $E_{ij}$ and $X_{i,:}h$ are both Gaussian random variables, independent of each other. Employing a union bound enables us to write for any $\tau > 0$
\begin{align*}
\mathbb{P}\Big[|m_i| \geq \tau\Big] &= \mathbb{P}\Big[|E_{ij}(X_{i,:}h)| \geq \tau\Big]\\
&\leq \mathbb{P}\Big[|E_{ij}| \geq \sqrt{\tau}\Big] + \mathbb{P}\Big[|X_{i,:}h| \geq \sqrt{\tau}\Big]\\
&\leq 2e^{-c_1\tau} + 2e^{-c_2\tau}\\
&\leq 4e^{-c\tau},
\end{align*}
for some suitable constants $c_1>0$, $c_2>0$ and $c>0$.

To establish that $\frac{1}{N}\sum_i m_i$ is sub-exponential, we note that the sub-exponential norm $\|\cdot\|_{\psi_1}$ (see \cite{vershynin10}, Definition 5.13) of $m_i$ can be upper-bounded by a constant.  We denote by $\kappa > 0$ the largest of these constants, i.e., 
\begin{align*}
\kappa = \underset{i=1,\ldots, N}{\max}\|m_i\|_{\psi_1} = \underset{i=1,\ldots, N}{\max}\|X_{i,:}h\|_{\psi_1}.
\end{align*}
%


Now, using  Lemma 5.15 in \cite{vershynin10}, the moment generating function of $m_i$ satisfies the following result: for $s$ such that $|s|\leq \frac{\eta}{\kappa}$ and for all $i = 1, \ldots, N$
\begin{align}
\label{eq:subexptail}
\mathbb{E}\Big[e^{sm_i}\Big] \leq e^{cs^2\kappa^2},
\end{align}
where $c$ and $\eta$ are absolute constants. Next, using Markov inequality, we can write for any $\varepsilon^{\prime} > 0$
\begin{align}
\label{eq:markov}
\mathbb{P}\left[\sum_{i=1}^N m_i \geq \varepsilon^\prime\right] &= \mathbb{P}\left[\exp\left(s\sum_{i=1}^N m_i\right) \geq \exp(s\varepsilon^\prime)\right]\nonumber\\
&\leq \frac{\mathbb{E}\left[\exp\left(s\sum_{i=1}^N m_i\right)\right]}{\exp(s\varepsilon^\prime)}.
\end{align}
To bound the numerator, we use \eqref{eq:subexptail} and write for $|s| \leq \frac{\eta}{\kappa}$ utilizing the iterated expectation
\begin{align*}
\mathbb{E}\left[\exp\left(s\sum_{i=1}^N m_i\right)\right] &= \mathbb{E}\left[\exp(sm_N)\exp\left(s\sum_{i=1}^{N-1} m_i\right)\right]\\
&= \mathbb{E}_{m_1,\ldots,m_{N-1}}\left[\mathbb{E}_{m_N|m_1,\ldots,m_{N-1}}\left[\exp(sm_N)\exp\left(s\sum_{i=1}^{N-1} m_i\right)\right]\right]\\
&= \mathbb{E}_{m_1,\ldots,m_{N-1}}\left[\mathbb{E}_{m_N|m_1,\ldots,m_{N-1}}\Big[\exp(sm_N)\Big]\exp\left(s\sum_{i=1}^{N-1} m_i\right)\right]\\
&\overset{\text{using}~ \eqref{eq:subexptail}}{\leq} \exp(cs^2\kappa^2)\mathbb{E}_{m_1,\ldots,m_{N-1}}\left[\exp\left(s\sum_{i=1}^{N-1} m_i\right)\right]\\
&\leq \exp(cs^2\kappa^2)\exp(cs^2\kappa^2)\mathbb{E}_{m_1,\ldots,m_{N-2}}\left[\exp\left(s\sum_{i=1}^{N-2} m_i\right)\right]\\
\vdots\\
&\leq \exp(Ncs^2\kappa^2)
\end{align*}
Substituting back to \eqref{eq:markov}, we get for $|s| \leq \frac{\eta}{\kappa}$
\begin{align}
\label{eq:markov2}
\mathbb{P}\left[\sum_{i=1}^N m_i \geq \varepsilon^\prime\right] \leq \exp(-s\varepsilon^\prime + Ncs^2\kappa^2).
\end{align}
We now select $s$ to minimize the right hand side of \eqref{eq:markov2}. For this, note that if the minimum is achieved for an $s$, which satisfies $|s| \leq \frac{\eta}{\kappa}$, then we simply minimize $-s\varepsilon^\prime + Ncs^2\kappa^2$ and get $s =\frac{\varepsilon^\prime}{N2c\kappa^2}$. On the other hand, if the minimum is achieved for an $s$ outside the range $|s| \leq \frac{\eta}{\kappa}$, we pick the one on boundary $s = \frac{\eta}{\kappa}$. Thus, choosing $s = \min\left(\frac{\varepsilon^\prime}{N2c\kappa^2}, \frac{\eta}{\kappa}\right)$, we obtain
\begin{align*}
\mathbb{P}\left[\sum_{i=1}^N m_i \geq \varepsilon^\prime\right] \leq \exp\left(-\min\left(\frac{{\varepsilon^\prime}^2}{4cN\kappa^2}, \frac{\eta\varepsilon^\prime}{2\kappa}\right)\right).
\end{align*}
Finally, setting $\varepsilon^\prime = N\varepsilon$, for a suitable constant $c > 0$, we get 
\begin{align*}
\mathbb{P}\left[\frac{1}{N}\sum_{i=1}^N m_i \geq \varepsilon\right] \leq \exp\left(-c \min\left(\frac{N\varepsilon^2}{\kappa^2}, \frac{N\varepsilon}{\kappa}\right)\right).
\end{align*}
Repeating the above argument for $- \frac{1}{N}\sum_{i=1}^N m_i$, we obtain same bound and a combination of both of them gives the required concentration inequality for the sum over the martingale difference sequence
\begin{align}
\label{eq:subexpbound}
\mathbb{P}\left[\frac{1}{N}\Bigg|\sum_{i=1}^N m_i\Bigg| \geq \varepsilon\right] = \mathbb{P}\left[\frac{1}{N}\Bigg|\left<X^TE_{:,j},h\right> \Bigg| \geq \varepsilon \right] \leq 2\exp\left(-c \min\left(\frac{N\varepsilon^2}{\kappa^2}, \frac{N\varepsilon}{\kappa}\right)\right).
\end{align}

\subsection{Establishing bound on $\mathbb{E}\Bigg[\underset{R(u_j)\leq r_j}{\sup}\frac{1}{N}\left<X^TE_{:,j},u_j\right>\Bigg]$}

To establish a high probability bound on the mean of $\underset{R(u_j) \leq r_j}{\text{sup}}\left<X^TE_{:,j},u_j\right>$, we use a generic chaining argument from \cite{talagrand06}, in particular Theorem 1.2.7 in \cite{tala05}. For this, we define $(Y_{u_j})_{u_j\in R(u_j)\leq r_j} = \frac{1}{N}\left<X^TE_{:,j},u_j\right>$ and $(Y_{v_j})_{v_j\in R(v_j)\leq r_j} = \frac{1}{N}\left<X^TE_{:,j},v_j\right>$ to be two centered random symmetric process, indexed by a fixed vectors $u_j$ and $v_j$, respectively. They are centered due to \eqref{eq:meanMi} and they are symmetric since, for example, the process $(Y_{u_j})_{u_j\in R(u_j)\leq r_j}$ has the same law as process $\Big(-(Y_{u_j})_{u_j\in R(u_j)\leq r_j}\Big)$ (see the results established in \eqref{eq:subexpbound}).  Consider now the absolute difference of these two processes
\begin{align*}
\Big|(Y_{u_j})_{u_j\in R(u_j)\leq r_j} - (Y_{v_j})_{v_j\in R(v_j)\leq r_j}\Big| &= \frac{1}{N}\Bigg|\left<X^TE_{:,j},u_j-v_j\right>\Bigg|\\
&= \|u_j-v_j\|_2\frac{1}{N}\Bigg|\left<X^TE_{:,j},\frac{u_j-v_j}{\|u_j-v_j\|_2}\right>\Bigg|.
\end{align*}
Using now the bound obtained in \eqref{eq:subexpbound}, we get
\begin{align*}
&\mathbb{P}\left[\frac{1}{N}\Bigg|\left<X^TE_{:,j},\frac{u_j-v_j}{\|u_j-v_j\|_2}\right> \Bigg| \geq \varepsilon \right]\\
= &\mathbb{P}\left[\|u_j-v_j\|_2\frac{1}{N}\Bigg|\left<X^TE_{:,j},\frac{u_j-v_j}{\|u_j-v_j\|_2}\right> \Bigg| \geq \|u_j-v_j\|_2\varepsilon \right]\\
= &\mathbb{P}\left[\frac{1}{N}\Bigg|\left<X^TE_{:,j}, u_j-v_j\right> \Bigg| \geq \tau \right]
 \leq 2\exp\left(-c \min\left(\frac{N\tau^2}{\|u_j-v_j\|_2^2\kappa^2}, \frac{N\tau}{\|u_j-v_j\|_2\kappa}\right)\right),
\end{align*}
where $\tau = \|u_j-v_j\|_2\varepsilon$. Then, according to Theorem 1.2.7 in \cite{tala05}, we obtain the following bound on the expectation of the supremum of the difference between the processes
\begin{align}
\label{eq:boundSupDiff}
\mathbb{E}\Bigg[\underset{R(u_j)\leq r_j, R(v_j)\leq r_j}{\sup}&\frac{1}{N}\Bigg|\left<X^TE_{:,j},u_j\right> - \left<X^TE_{:,j},v_j\right>\Bigg|\Bigg]\nonumber\\
 &\hspace{20pt}\leq c\left(\gamma_1\left(S_j, \frac{\|u_j-v_j\|_2}{N}\right) + \gamma_2\left(S_j, \frac{\|u_j-v_j\|_2}{\sqrt{N}}\right)\right),
\end{align}
where $c$ is a constant, $f_i(S_j, d_i)$, $i=1,2$, are the majorizing measures, which are defined in \cite{talagrand06}, Definition 1.2.5; $d_1 = \frac{\|u_j-v_j\|_2}{N}$ and  $d_2 = \frac{\|u_j-v_j\|_2}{\sqrt{N}}$ are the distance measures on the set $S_j$ defined for all vectors $s \in S_j: R(s)\leq r_j$. The definition of majorizing measure is as follows, for $\alpha > 0$
\begin{align}
\label{eq:gammaFun}
\gamma_\alpha(S_j, d) = \inf \underset{t}{\sup}\sum_{k\geq 0} 2^{\frac{k}{\alpha}}\Delta(A_k(t)),
\end{align}
where $\inf$ is taken over all possible admissible sequences of the set $S_j$; $\Delta(A_k(t))$ denotes the diameter of element $A_k(t)$ with respect to the distance metric $d$ defined as
\begin{align}
\label{eq:diameter}
\Delta(A_k(t)) = \underset{t_1, t_2 \in A_k(t)}{\sup}d(t_1,t_2),
\end{align} 
and $A_k(t) \in \mathcal{A}_k$ is an element of an admissible sequence in generic chaining, see Definition 1.2.3 in \cite{talagrand06} for a detailed discussion on how $\mathcal{A}_k$ are constructed.

Observe that from definition of a diameter $\Delta(\cdot)$ in \eqref{eq:diameter} and majorizing measure in \eqref{eq:gammaFun} we can immediately see that for any constant $c > 0$
\begin{align}
\label{eq:f1f1ineq}
\gamma_\alpha\left(S_j, cd\right) = c\gamma_\alpha\left(S_j, d\right),
\end{align}
since $\inf \underset{t}{\sup}\sum_{k\geq 0} 2^{\frac{k}{\alpha}}\underset{t_1, t_2 \in A_k(t)}{\sup}cd(t_1,t_2) = c\inf \underset{t}{\sup}\sum_{k\geq 0} 2^{\frac{k}{\alpha}}\underset{t_1, t_2 \in A_k(t)}{\sup}d(t_1,t_2)$.  Moreover, in the next result we establish the following useful Lemma which would enable us to bound the $\gamma_1$ with the square of $\gamma_2$. 

\begin{lemma}
\label{eq:lemMaj}
Given a metric space $(S_j,d)$, we have
\begin{align}
\gamma_1(S_j, \|.\|_2) \leq \gamma_2^2(S_j, \|.\|_2). 
\end{align}
\label{}
 \end{lemma}

To prove this Lemma, we define $d(s,t) = \|s-t\|_2$. We use the traditional definition of majorizing measure $\gamma_{\alpha}^\prime(S_j,d)$ from \cite{tala01}, equation (1.2): 
\begin{align*}
\gamma_{\alpha}^\prime(S_J,d) = \inf \underset{s\in S}{\sup} \left( \int_0^\infty \left( \log \frac{1}{\mu(B_d(s,\varepsilon))} \right)^{1/\alpha} d \varepsilon \right),
\end{align*}
where $B_d(s,\varepsilon)$ is the closed ball of center $t$ and radius $\varepsilon$ based on the distance $d$ and the infimum is taken over all the probability measure $\mu$ on $S_j$.
 
Note that $\gamma_{\alpha}^\prime(S_j,d)$ relates to the majorizing measure $\gamma_{\alpha}(S_j,d)$ used in \eqref{eq:boundSupDiff}  as (see \cite{tala01}, Theorem 1.2)
\begin{align*}
 K(\alpha)^{-1} \gamma_\alpha(S_j,d) \leq \gamma_{\alpha}^\prime(S_j,d) \leq K(\alpha) \gamma_\alpha(S_j,d),
\end{align*}
where $K(\alpha)$ is a constant depending on $\alpha$ only. As a result, it is enough to show that $\gamma_{1}^\prime(S_j, d) \leq {\gamma_{2}^\prime}^2(S_j, d)$. The required relationship is then established as follows
\begin{align*}
\gamma_{1}^\prime(S_j,d) &= \inf \sup_t \left( \int_0^\infty \left( \log \frac{1}{\mu(B_{d}(t,\varepsilon))} \right) d \varepsilon \right) \\
& \leq \inf \sup_t \left( \int_0^\infty \left( \log \frac{1}{\mu(B_{d}(t,\varepsilon))} \right)^{1/2} d \varepsilon \right)^2 \\
&= {\gamma_{2}^\prime}^2(S_j,d).
\end{align*}
And this completes the proof. Now using Theorem 2.1.1 in \cite{talagrand06}, and the definition of $\gamma_\alpha(S_j,d)$ in \eqref{eq:gammaFun} we can establish that
\begin{align}
\label{eq:f2}
\gamma_2\left(S_j, \frac{\|.\|_2}{\sqrt{N}}\right) &= \frac{1}{\sqrt{N}}\gamma_2(S_j, \|.\|_2)\nonumber \quad\text{using \eqref{eq:f1f1ineq}}\\
 &\leq \frac{1}{\sqrt{N}}\mathbb{E}\Big[\underset{R(z) \leq r_j}{\sup}\left<g,z\right>\Big]\nonumber \quad\text{using Theorem 2.1.1 in \cite{talagrand06}}\\
&=r_j\frac{1}{\sqrt{N}}\mathbb{E}\Big[\underset{R(u) \leq 1}{\sup}\left<g,u\right>\Big]\nonumber\quad\text{since $\mathbb{E}\Big[\underset{R(z) \leq r_j}{\sup}\left<g,z\right>\Big] = r_j\mathbb{E}\Big[\underset{R(u) \leq 1}{\sup}\left<g,u\right>\Big]$ for $z=r_ju$}\\
&=r_j\frac{1}{\sqrt{N}}w(\Omega_R),
\end{align}
where in the last line we used the description of Gaussian width in Definition \ref{eq:gaussWidth}. Using Lemma \ref{eq:lemMaj} and \eqref{eq:f1f1ineq} above, we also get
\begin{align}
\label{eq:f1}
\gamma_1\left(S_j, \frac{\|.\|_2}{N}\right) &= \frac{1}{N}\gamma_1\left(S_j, \|.\|_2\right)\nonumber\quad\text{using \eqref{eq:f1f1ineq}}\\
 &\leq \frac{1}{N}\gamma_2^2\left(S, \|.\|_2\right)\nonumber\quad\text{using Lemma \ref{eq:lemMaj}}\\
&\leq r_j^2\frac{1}{N^2}w^2(\Omega_R)\nonumber\quad\text{using \eqref{eq:f2}}\\
&\leq r_j\frac{1}{N^2}w^2(\Omega_R),
\end{align}
where in the last line we used the fact that $r_j < 1$. Finally, substituting \eqref{eq:f2} and \eqref{eq:f1} into \eqref{eq:boundSupDiff} and using Lemma 1.2.8 in \cite{tala05}, we get
\begin{align}
\label{eq:boundESup}
\mathbb{E}\Bigg[\underset{R(u_j)\leq r_j, R(v_j)\leq r_j}{\sup}\frac{1}{N}\Bigg|\left<X^TE_{:,j},u_j\right> - \left<X^TE_{:,j},v_j\right>\Bigg|\Bigg] = \mathbb{E}\Bigg[\underset{R(u_j)\leq r_j}{\sup}\Bigg|\frac{1}{N}\left<X^TE_{:,j},u_j\right> \Bigg|\Bigg]\nonumber\\
 \leq cr_j\left(\frac{w(\Omega_R)}{\sqrt{N}} + \frac{w^2(\Omega_R)}{N^2}\right).
\end{align}

\subsection{Establishing high probability concentration bound}
Next, in order to establish the high probability concentration of  the supremum of the random variable $\frac{1}{N}\left<X^TE_{:,j},u_j\right>$ around its mean, we use Theorem 1.2.9 from \cite{talagrand06}. For any $\epsilon_1 > 0$ and $\epsilon_2 > 0$, we have
\begin{align}
\label{eq:concentrationLambda}
\mathbb{P}\left[\underset{R(u_j)\leq r_j}{\sup}\Bigg|\frac{1}{N}\left<X^TE_{:,j},u_j\right>\Bigg| \geq \mathbb{E}\Bigg[\underset{R(u_j)\leq r_j}{\sup}\Bigg|\frac{1}{N}\left<X^TE_{:,j},u_j\right> \Bigg|\Bigg] + \epsilon_1 D_1 + \epsilon_2 D_2\right] \leq c \exp(-\min(\epsilon_2^2, \epsilon_1)).
\end{align}
where $D_i \leq \gamma_i(S_j, d)$, $i=1,2$, where $\gamma_i(S_j, d)$ are as defined in the discussion after \eqref{eq:boundSupDiff}. Therefore, using the result \eqref{eq:boundESup}, the concentration inequality \eqref{eq:concentrationLambda} can now be written as 
\begin{align}
\label{eq:concentrationLambda2}
\mathbb{P}\Bigg[\underset{R(u_j) \leq r_j}{\text{sup}}\Bigg|\frac{1}{N}\left<X^TE_{:,j},u_j\right>\Bigg| \geq \left(c_2(1+\epsilon_2)r_j\frac{w(\Omega_R)}{\sqrt{N}} + c_1(1+\epsilon_1)r_j\frac{w^2(\Omega_R)}{N^2}\right) \Bigg] \leq c \exp(-\min(\epsilon_2^2, \epsilon_1)).
\end{align}
To adapt to the form required in \eqref{eq:objectiveInequlityy}, we reverse the direction of inequality
 \begin{align}
 \label{eq:concentrationLambda3}
 \mathbb{P}\Bigg[\underset{R(u_j) \leq r_j}{\text{sup}}\Bigg|\frac{1}{N}\left<X^TE_{:,j},u_j\right>\Bigg| \leq \left(c_2(1+\epsilon_2)r_j\frac{w(\Omega_R)}{\sqrt{N}} + c_1(1+\epsilon_1)r_j\frac{w^2(\Omega_R)}{N^2}\right) \Bigg]\nonumber\\
  \geq 1- c \exp(-\min(\epsilon_2^2, \epsilon_1)).
 \end{align}

\subsection{Overall bound}
\label{sec:lambdaOverallBound}

Now we can combine the results obtained in \eqref{eq:concentrationLambda3} for each $j=1,\ldots, p$ using the fact that $\sum_{j=1}^pr_j \leq 1$ and using the form of the overall bound in \eqref{eq:multiStepBound}. Therefore, we get
\begin{align*}
\mathbb{P}\Bigg[\underset{R(U) \leq 1}{\text{sup}} \left<\frac{1}{N}Z^T\bm{\epsilon}, U\right> \leq \left(c_2(1+\epsilon_2)\frac{w(\Omega_R)}{\sqrt{N}} + c_1(1+\epsilon_1)\frac{w^2(\Omega_R)}{N^2}\right) \Bigg]\\ \geq 1- c \exp(-\min(\epsilon_2^2, \epsilon_1) + \log(p)).
\end{align*}
This concludes our proof on establishing the bound on the regularization parameter.

\section{Restricted Eigenvalue Condition}
\label{sec:re_append}
To establish restricted eigenvalue (RE) condition, we need to show that $\frac{||(I_{p\times p}\otimes X)\Delta||_2}{||\Delta||_2} \geq \sqrt{\kappa N}$, $\kappa > 0$, for all $\Delta = \hat{\bm{\beta}} - \bm{\beta}^*$, $\Delta \in \text{cone}(\Omega_E)$, where $\text{cone}(\Omega_E)$ denotes a cone of an error set\\ $\Omega_E = \left\{ \Delta \in \mathbb{R}^{dp^2} \Big| R(\bm{\beta}^* + \Delta) \leq R(\bm{\beta}^*) + \frac{1}{c}R(\Delta)\right\}$. To show $\frac{||(I_{p\times p}\otimes X)\Delta||_2}{||\Delta||_2} \geq \sqrt{\kappa N}$ for all $\Delta \in \text{cone}(\Omega_E)$, we will show that $\underset{\Delta \in \text{cone}(\Omega_E)}{\inf}\frac{||(I_{p\times p}\otimes X)\Delta||_2}{||\Delta||_2} \geq \sqrt{\rho}$, for some $\rho > 0$ and then set $\kappa N = \rho$.

Note that the error vector can be written as $\Delta = [\Delta_1^T, \Delta_2^T, \ldots, \Delta_p^T]^T$, where $\Delta_i$ is of size $dp\times 1$. Also let $\bm{\beta}^* = [\beta_1^{*T} \beta_2^{*T} \ldots \beta_p^{*T}]^T$, for $\beta_i^* \in \mathbb{R}^{dp}$, then using our assumption in \eqref{eq:decompNorm} that the norm $R(\cdot)$ is decomposable, we can represent original set $\Omega_E$ as a Cartesian product of subsets $\Omega_{E_i}$, i.e., $\Omega_E = \Omega_{E_1} \times \Omega_{E_2} \times \cdots \times \Omega_{E_p}$, where
\begin{align*}
\Omega_{E_i} = \left\{\Delta_i \in \mathbb{R}^{dp} \Big| R(\beta_i^* + \Delta_i) \leq R(\beta_i^*) + \frac{1}{c}R(\Delta_i) \right\},
\end{align*}
which also implies that $\text{cone}(\Omega_E) = \text{cone}(\Omega_{E_1})\times\text{cone}(\Omega_{E_2})\times \cdots \times \text{cone}(\Omega_{E_p})$.
Also, if $||\Delta||_2 = 1$, then we denote $||\Delta_i||_2 = \delta_i > 0$, so that $\sum_{i=1}^p\delta_i^2 = 1$. With this information, we can write
\begin{align}
\label{eq:REIneq}
\underset{\Delta \in \text{cone}(\Omega_E)}{\inf}\frac{||(I_{p\times p}\otimes X)\Delta||_2^2}{||\Delta||_2^2} 
&= \underset{\begin{subarray}{c}
\Delta \in \text{cone}(\Omega_E)\\
||\Delta||_2 = 1
\end{subarray}
}{\inf}||(I_{p\times p}\otimes X)\Delta||_2^2\nonumber\\
&= \underset{\begin{subarray}{c}
 \Delta \in \text{cone}(\Omega_E)\nonumber\\
 ||\Delta||_2 = 1
 \end{subarray}
 }{\inf}
 ||X\Delta_1||_2^2 + ||X\Delta_2||_2^2 + \ldots + ||X\Delta_p||_2^2 \\
 &= \sum_{i=1}^p~\underset{\begin{subarray}{c}
 \Delta_i \in \text{cone}(\Omega_{e_i})\\
 ||\Delta_i||_2 = \delta_i
 \end{subarray}
 }{\inf}||X\Delta_i||_2^2.
\end{align}
Our objective is to establish a high probability bound of the form
\begin{align*}
\mathbb{P}\Bigg[\underset{\Delta \in \text{cone}(\Omega_E)}{\inf}\frac{||(I_{p\times p}\otimes X)\Delta||_2}{||\Delta||_2} \geq \rho\Bigg] \geq \pi
\end{align*} 
where $0\leq \pi\leq 1$, i.e., lower bound should hold with at least probability $\pi$. Note that if we square the terms inside the probability statement above, the probability of the resulting expression does not change since the squared terms are positive. Therefore, using \eqref{eq:REIneq} and assuming that $\rho^2 = \sum_{i=1}^p\rho_i^2$ we can rewrite the above as follows
\begin{align}
\label{eq:multiStepBoundRE2}
\mathbb{P}\Bigg[\underset{\Delta \in \text{cone}(\Omega_E)}{\inf}\frac{||(I_{p\times p}\otimes X)\Delta||_2}{||\Delta||_2} \geq \rho\Bigg]
&=\mathbb{P}\Bigg[\underset{\Delta \in \text{cone}(\Omega_E)}{\inf}\frac{||(I_{p\times p}\otimes X)\Delta||_2^2}{||\Delta||_2^2} \geq \sum_{i=1}^p\rho_i^2\Bigg] \nonumber\\
&=
\mathbb{P}\Bigg[
\sum_{i=1}^p~\underset{\begin{subarray}{c}
 \Delta_i \in \text{cone}(\Omega_{E_i})\\
 ||\Delta_i||_2 = \delta_i
 \end{subarray}
 }{\inf}||X\Delta_i||_2^2\geq \sum_{i=1}^p\rho_i^2\Bigg]\nonumber\quad\text{using \eqref{eq:REIneq}}\\
&\geq
\mathbb{P}\Bigg[ \Bigg\{\underset{\begin{subarray}{c}
 \Delta_1 \in \text{cone}(\Omega_{E_1})\\
 ||\Delta_1||_2 = \delta_1
 \end{subarray}
 }{\inf}||X\Delta_1||_2^2\geq \rho_i^2\Bigg\} \text{ and }\nonumber \\
  &\hspace{60pt}\ldots \text{ and }
 \Bigg\{\underset{\begin{subarray}{c}
  \Delta_p \in \text{cone}(\Omega_{E_p})\\
  ||\Delta_p||_2 = \delta_p
  \end{subarray}
  }{\inf}||X\Delta_p||_2^2\geq \rho_i^2\Bigg\}
  \Bigg] \nonumber\\
&\geq
\sum_{i=1}^p
\mathbb{P}\Bigg[ \underset{\begin{subarray}{c}
 \Delta_i \in \text{cone}(\Omega_{E_i})\\
 ||\Delta_i||_2 = \delta_i
 \end{subarray}
 }{\inf}||X\Delta_i||_2^2\geq \rho^2_i\Bigg] - (p-1)\nonumber\\
&=
\sum_{i=1}^p
\mathbb{P}\Bigg[ \underset{\begin{subarray}{c}
 \Delta_i \in \text{cone}(\Omega_{E_i})\\
 ||\Delta_i||_2 = \delta_i
 \end{subarray}
 }{\inf}||X\Delta_i||_2\geq \rho_i\Bigg] - (p-1)\nonumber \quad\text{taking square root}\\
 &=
 \sum_{i=1}^p 
  \mathbb{P}\Bigg[ \underset{\begin{subarray}{c}
   \Delta_i \in \text{cone}(\Omega_{E_i})\\
   ||\Delta_i||_2 = \delta_i
   \end{subarray}
   }{\inf}\frac{||X\Delta_i||_2}{||\Delta_i||_2}\geq \frac{\rho_i}{||\Delta_i||_2}\Bigg] - (p-1)\nonumber\\
 &= 
\sum_{i=1}^p \mathbb{P}\Bigg[ \underset{u_i \in \text{cone}(\Omega_{E_i})\cap S^{dp-1}}{\inf}||Xu_i||_2\geq \frac{\rho_i}{\delta_i}\Bigg] - (p-1)
\end{align}
where we defined $u_i = \frac{\Delta_i}{||\Delta_i||_2}$ and $S^{dp-1}$ is a unit sphere.
Therefore, if we denote $\Theta_i  = \text{cone}(\Omega_{E_i})\cap S^{dp-1}$, we need to establish a lower bound of the form
\begin{align}
\label{eq:REboundTypee}
\mathbb{P}\Bigg[ \underset{u_i \in \Theta_i}{\inf}||Xu_i||_2\geq \rho_i^\prime\Bigg] \geq \pi_i,
\end{align}
where $\rho_i^\prime = \frac{\rho_i}{\delta_i}$. In the following derivations we set $\Theta = \text{cone}(\Omega_{E_i})\cap S^{dp-1}$ and  $u = u_i$ for all $i=1,\ldots, p$ since the specific index $i$ is irrelevant.

\subsection{Bound on $\underset{u \in \Theta}{\inf }~ ||Xu||_2$}

Using results from Appendix \ref{sec:covariance_append} we can establish that $Xu \in \mathbb{R}^{N}$ is a Gaussian random vector, i.e., $Xu \sim \mathcal{N}(0, Q_u)$, where covariance matrix $Q_u = (I_{N\times N}\otimes u^T)C_\mathcal{U}(I_{N\times N}\otimes u)$, $C_\mathcal{U}$ is defined in \eqref{eq:covU}, and  $u \in \Theta$ is a fixed vector. 

To establish $\underset{u \in \Theta}{\inf }~ ||Xu||_2$, we invoke a generic chaining argument from \cite{talagrand06}, specifically Theorem 2.1.5.  For this we let $(Z_u)_{u \in \Theta} = ||Xu||_2 - \mathbb{E}(||Xu||_2)$ and $(Z_v)_{v \in \Theta} = ||Xv||_2 - \mathbb{E}(||Xv||_2)$ be two centered symmetric random processes. They are centered since, for example, $\mathbb{E}\Big[(Z_u)_{u \in \Theta}\Big] = \mathbb{E}(||Xu||_2) - \mathbb{E}(||Xu||_2) = 0$, and they are symmetric due to the later result shown in \eqref{eq:subGaussian}. 

\vspace{15pt}
\noindent\textbf{Sub-gaussianity of the process $Z_u - Z_v$}.\newline We can show that the process difference
\begin{align}
\label{eq:processDiff}
(Z_u)_{u \in \Theta} - (Z_v)_{v \in \Theta} = \|u-v\|_2\left(\left\|X\frac{u-v}{\|u-v\|_2}\right\|_2 - \mathbb{E}\left(\left\|X\frac{u-v}{\|u-v\|_2}\right\|_2\right)\right)
\end{align}
is a sub-Gaussian random process. This is indeed the case since we can establish that for $Z = ||X\frac{u-v}{\|u-v\|_2}||_2 - \mathbb{E}(||X\frac{u-v}{\|u-v\|_2}||_2)$, the sub-gaussian norm  $\|Z\|_{\psi_2}\leq K$ for some constant $K > 0$ (see \cite{vershynin10}, Definition 5.7). To show this, let $\xi  = \frac{u-v}{\|u-v\|_2}$ and apply concentration of a Lipschitz function of Gaussian random variables. Specifically, observe that $X\xi \sim \mathcal{N}(0, Q_\xi)$ is distributed same as  $\sqrt{Q_\xi}g \sim \mathcal{N}(0, Q_\xi)$, where $g\sim \mathcal{N}(0, I_{N\times N})$. Therefore, we can write 
\begin{align*}
\mathbb{P}\Big[\left|\|X\xi\|_2 - \mathbb{E}(\|X\xi\|_2)\right|>\tau\Big]  = \mathbb{P}\Big[\left|\|\sqrt{Q_\xi}g\|_2 - \mathbb{E}(\|\sqrt{Q_\xi}g\|_2)\right|>\tau\Big].
\end{align*}
Moreover, note that $\|\sqrt{Q_\xi}g\|_2$ is a Lipschitz function with constant $\|\sqrt{Q_\xi}\|_2$ since we can write $\Big|\|\sqrt{Q_\xi}g_1\|_2 - \|\sqrt{Q_\xi}g_2\|_2\Big| \leq \|\sqrt{Q_\xi}(g1 -g2)\|_2 \leq \|\sqrt{Q_\xi}\|_2~\|g_1-g_2\|_2$. Using the concentration of a Lipschitz function of Gaussian random variables, we can obtain for all $\tau > 0$
\begin{align}
\label{eq:subGaussian}
\mathbb{P}\Big[\left|\|X\xi\|_2 - \mathbb{E}(\|X\xi\|_2)\right|>\tau\Big]  &= \mathbb{P}\Big[\left|\|\sqrt{Q_\xi}g\|_2 - \mathbb{E}(\|\sqrt{Q_\xi}g\|_2)\right|>\tau\Big]\nonumber\\
&\leq 2\exp\left(-\frac{\tau^2}{2\|Q_{\xi}\|_2}\right)\nonumber\\
&\leq 2\exp\left(-\frac{\tau^2}{2\mathscr{M}}\right),
\end{align}
where $||Q_\xi||_2 \leq ||\xi||_2^2\frac{\Lambda_{\max}(\Sigma)}{\Lambda_{\min}(\bm{\mathscr{A}})} = \frac{\Lambda_{\max}(\Sigma)}{\Lambda_{\min}(\bm{\mathscr A})} = \mathscr{M}$ (see \eqref{eq:Qnorm2}), and which shows that $\|X\xi\|_2$ is sub-Gaussian with constant $K = \sqrt{\mathscr{M}}$. 

Now, using \eqref{eq:subGaussian} we can establish the sub-Gaussian tails of \eqref{eq:processDiff}. Define $\tau^\prime =\|u-v\|_2\tau$ and write 
\begin{align}
\label{eq:subGaussian2}
\mathbb{P}\Big[\left| \|u-v\|_2\Big(\|X\xi\|_2 - \mathbb{E}(\|X\xi\|_2)\Big)\right|> \|u-v\|_2\tau\Big] &= \mathbb{P}\Big[\left|(Z_u)_{u \in \Theta} - (Z_v)_{v \in \Theta}\right|>\tau^\prime\Big] \nonumber\\
&\leq 2\exp\left(-\frac{{\tau^\prime}^2}{2\|u-v\|_2^2\mathscr{M}}\right).
\end{align}

\noindent \textbf{Establishing bound on $\mathbb{E}\left(\underset{u \in \Theta}{\inf}||Xu||_2\right)$}.\newline 
Using the results established in \eqref{eq:subGaussian2} and Theorem 2.1.5 in \cite{talagrand06}, we can conclude that the distance measure on the set $\Theta$ is $d(u,v) = \|u-v\|_2$ for $u,v \in \Theta$. Moreover, we can now obtain an upper bound on the expectation of the supremum of the process difference $|Z_u-Z_v|$
\begin{align}
\label{eq:REsup}
\mathbb{E}\left(\underset{u,v \in \Theta}{\sup}\Big|Z_u - Z_v\Big|\right) &= \mathbb{E}\left(\underset{u,v \in \Theta}{\sup}\Big| ~||X(u-v)||_2 - \mathbb{E}(||X(u-v)||_2)~\Big|\right)\nonumber\\
&= \mathbb{E}\left(\underset{u \in \Theta}{\sup}\Big| ~||Xu||_2 - \mathbb{E}(||Xu||_2)~\Big|\right)\nonumber\quad\text{using Lemma 1.2.8 in \cite{talagrand06}}\\
&\leq \mathbb{E}\Big[\underset{u \in \Theta}{\sup}\left<g,u\right>\Big]\nonumber\\
&\leq c w(\Theta),
\end{align}
where $g \sim \mathcal{N}(0, I)$, $w(\Theta)$ is the Gaussian width of set $\Theta$ and $c$ is a constant.

Since we are interested in the bound on $\underset{u \in \Theta}{\inf }~ ||Xu||_2$, we can extract from \eqref{eq:REsup} the lower bound on the expectation of the infimum of the process. Specifically, note that \eqref{eq:REsup} can be written as  
\begin{align*}
\mathbb{E}\left(\left| \underset{u \in \Theta}{\inf}~||Xu||_2 - \underset{u \in \Theta}{\inf}\mathbb{E}(||Xu||_2)~\right|\right) \leq \mathbb{E}\left(\underset{u \in \Theta}{\sup}\bigg| ~||Xu||_2 - \mathbb{E}(||Xu||_2)~\bigg|\right) \leq c w(\Theta),
\end{align*}
leading to 
\begin{align*}
-c w(\Theta) \leq \mathbb{E}\left(\underset{u \in \Theta}{\inf}~||Xu||_2 - \underset{u \in \Theta}{\inf}\mathbb{E}(||Xu||_2)\right) \leq c w(\Theta).
\end{align*}
The lower bound then takes the form
\begin{align}
\label{eq:REinf}
\mathbb{E}\left(\underset{u \in \Theta}{\inf}||Xu||_2\right) \geq \underset{u \in \Theta}{\inf}~\mathbb{E}(||Xu||_2) - c w(\Theta)
\end{align}
Note that the vector $Xu$ is distributed  as $Xu \sim \mathcal{N}(0, Q_u)$, which is the same as a vector $\sqrt{Q_u}g \sim \mathcal{N}(0, Q_u)$ for $g \sim \mathcal{N}(0, I)$. Therefore, using results of Lemma I.2 from \cite{negahban11}, we can extract the following inequality
\begin{align*}
\Big|\sqrt{\text{trace}(Q_u)} - \mathbb{E}(\|\sqrt{Q_u}g\|_2)\Big| \leq  2\sqrt{\Lambda_{\max}(Q_u)}.
\end{align*}
Moreover, based on our discussion, the same inequality holds for the random vector $Xu$ since $\mathbb{E}(\|\sqrt{Q_u}g\|_2) =  \mathbb{E}(\|Xu\|_2)$
\begin{align*}
\Big|\sqrt{\text{trace}(Q_u)} - \mathbb{E}(\|Xu\|_2)\Big| \leq  2\sqrt{\Lambda_{\max}(Q_u)}.
\end{align*}
which leads to a lower bound on the expectation of the norm 
\begin{align}
\label{eq:EXu}
\mathbb{E}(\|Xu\|_2) \geq  \sqrt{\text{trace}(Q_u)} - 2\sqrt{\Lambda_{\max}(Q_u)}.
\end{align}
We will lower-bound the first term on the right hand side of \eqref{eq:EXu} and upper bound the second one. In particular, using \eqref{eq:Qtrace} we write ${\text{trace}}(Q_u) = Nu^TC_\mathsf{X}u$ for any $u\in \Theta$ and bound
\begin{align}
\label{eq:bound2}
{\text{trace}}(Q_u) = Nu^TC_\mathsf{X}u &= N||C_\mathsf{X}^{\frac{1}{2}}u||_2^2 \nonumber \\
 &\geq N\underset{u \in \Theta}{\inf}~u^TC_\mathsf{X}u\nonumber \\
 &\geq N\underset{u \in \mathbb{R}^{dp}}{\inf}~u^TC_\mathsf{X}u  = N\Lambda_{\min}(C_\mathsf{X})\nonumber \\
 &\geq N\frac{\Lambda_{\min}(\Sigma)}{\Lambda_{\max}(\mathscr{A})} = N\mathscr{L}. 
\end{align}
Moreover, using \eqref{eq:Qnorm2}, we bound 
\begin{align}
\label{eq:bound1}
||Q_u||_2 \leq ||u||_2^2\frac{\Lambda_{\max}(\Sigma)}{\Lambda_{\min}(\bm{\mathscr{A}})} = \frac{\Lambda_{\max}(\Sigma)}{\Lambda_{\min}(\bm{\mathscr A})} = \mathscr{M}.
\end{align} 
Therefore, substituting \eqref{eq:bound1} and \eqref{eq:bound2} into \eqref{eq:EXu}, we get
\begin{align*}
\mathbb{E}(\|Xu\|_2) \geq \sqrt{N\mathscr{L}} - 2\sqrt{\mathscr{M}}.
\end{align*}
Since $\mathbb{E}(\|Xu\|_2)$ is bounded from below, we can write
\begin{align}
\label{eq:Exu2}
\underset{u \in \Theta}{\inf}~\mathbb{E}(\|Xu\|_2) \geq \sqrt{N\mathscr{L}} - 2\sqrt{\mathscr{M}}.
\end{align}
Finally, substituting \eqref{eq:Exu2} in \eqref{eq:REinf} gives us
\begin{align}
\label{eq:REinf2}
\mathbb{E}\left(\underset{u \in \Theta}{\inf}||Xu||_2\right) \geq \sqrt{N\mathscr{L}} - 2\sqrt{\mathscr{M}} - cw(\Theta).
\end{align}

\vspace{15pt}
\noindent \textbf{Establishing concentration inequality of $\underset{u \in \Theta}{\inf}||Xu||_2$}.\\
Now from Lemma 2.1.3 in \cite{talagrand06} and the results in \cite{banerjee14} we extract the form of the high probability concentration inequality of $\underset{u \in \Theta}{\inf }~ ||Xu||_2$ around its mean, for $\tau > 0$
\begin{align*}
\mathbb{P}\left[\underset{u \in \Theta}{\inf}||Xu||_2 \leq \mathbb{E}\left(\underset{u \in \Theta}{\inf}||Xu||_2\right) - \tau\right] \leq c_1\exp(-c_2 \tau^2).
\end{align*} 
In order to bring the above expression into the form of \eqref{eq:REboundTypee}, we write 
\begin{align*}
\mathbb{P}\left[\underset{u \in \Theta}{\inf}||Xu||_2 \geq \mathbb{E}\left(\underset{u \in \Theta}{\inf}||Xu||_2\right) - \tau\right] \geq 1 - c_1\exp(-c_2 \tau^2).
\end{align*}
Substituting the bound on the expectation from \eqref{eq:REinf2} gives us
\begin{align}
\label{eq:REbound1}
\mathbb{P}\left[\underset{u \in \Theta}{\inf}||Xu||_2 \geq  \sqrt{N\mathscr{L}} - 2\sqrt{\mathscr{M}} - cw(\Theta) - \tau \right] \leq c_1\exp(-c_2\tau^2).
\end{align}

\subsection{Overall bound}

Observe that in \eqref{eq:REbound1} we established a bound for each $u_i = \frac{\Delta_i}{||\Delta_i||_2}$ of the form 
\begin{align*}
\mathbb{P}\left[
\underset{\begin{subarray}{c}
  \Delta_i \in \text{cone}(\Omega_{E_i})\\
  ||\Delta_i||_2 = \delta_i
  \end{subarray}
  }{\inf}
\frac{||X\Delta_i||_2}{||\Delta_i||_2} \geq ||\Delta_i||_2~\rho_i^\prime\right]
\geq 1 - c_1\exp(-c_2\eta_i^2),
\end{align*}
where $\rho_i^\prime = \sqrt{N\mathscr{L}} - 2\sqrt{\mathscr{M}} - c w(\Theta) - \eta_i$. Then using the fact that $\rho_i = \rho_i^\prime\delta_i$, ~$\rho^2 = \sum_{i=1}^p\rho_i^2$, ~$\sum_i^p\delta_i^2=1$ and setting $\eta_i = \eta$ for all $i=1,\ldots, p$, we get
\begin{align*}
\rho^2 = \bigg[\sqrt{N\mathscr{L}} - 2\sqrt{\mathscr{M}} - c w(\Theta) - \eta\bigg]^2\sum_{i=1}^p\delta_i^2 = \bigg[\sqrt{N\mathscr{L}} - 2\sqrt{\mathscr{M}} - c w(\Theta) - \eta\bigg]^2.
\end{align*}
Taking the square root of the above and using \eqref{eq:multiStepBoundRE2} we finally get
\begin{align}
\label{eq:REresult}
\mathbb{P}\Bigg[\underset{\Delta \in \text{cone}(\Omega_E)}{\inf}\frac{||(I_{p\times p}\otimes X)\Delta||_2}{||\Delta||_2} \geq 
\sqrt{N\mathscr{L}} - 2\sqrt{\mathscr{M}} - c w(\Theta) - \eta
\Bigg] \geq  1 - pc_1\exp(-c_2\eta^2).
\end{align}

\vspace{15pt}
\noindent \textbf{Establishing bound on $N$}.\\
Now setting $\eta = \varepsilon\sqrt{N\mathscr{L}}$ for $0<\varepsilon < 1$, the right hand side of the inequality inside the probability statement in \eqref{eq:REresult} must be equal to
\begin{align*}
\sqrt{\kappa N} = \sqrt{N\mathscr{L}} - 2\sqrt{\mathscr{M}} - cw(\Theta) - \varepsilon\sqrt{N\mathscr{L}} = \varepsilon^\prime\sqrt{N\mathscr{L}} - 2\sqrt{\mathscr{M}} - cw(\Theta),
\end{align*}
for some positive constant  $\varepsilon^\prime$.
Since $\kappa N > 0$, it follows that we require
\begin{align*}
\varepsilon^\prime\sqrt{N\mathscr{L}} > 2\sqrt{\mathscr{M}} + c w(\Theta),
\end{align*}
or equivalently
\begin{align*}
\sqrt{N} > \frac{2\sqrt{\mathscr{M}} + c w(\Theta)}{\varepsilon^\prime\sqrt{\mathscr{L}}} = \mathcal{O}(w(\Theta)).
\end{align*}
This concludes our proof on establishing the restricted eigenvalue conditions.

\vspace*{3mm}

{\bf Acknowledgements:} The research was supported by NSF grants IIS-1447566, IIS-1422557, CCF-1451986, CNS-
1314560, IIS-0953274, IIS-1029711, and by NASA grant NNX12AQ39A.

\bibliography{main}
\bibliographystyle{plain}

\end{document}